\documentclass{amsart}

\usepackage{amssymb}
\usepackage{verbatim,latexsym}
\usepackage{amscd}
\usepackage{amssymb}
\usepackage{amsmath}

\theoremstyle{plain}
\newtheorem{Lem}{Lemma}[section]
\newtheorem{Cor}[Lem]{Corollary}
\newtheorem{Thm}[Lem]{Theorem}
\theoremstyle{definition} 

\newtheorem{Rk}[Lem]{Remark}
\newtheorem{Def}[Lem]{Definition}

\newcommand{\GG}{G_n} 
\newcommand{\lbr}{[\negthinspace[}
\newcommand{\rbr}{]\negthinspace]}
\newcommand{\colim}{\mathrm{colim}} 
\newcommand{\zig}{\addtocounter{Lem}{1}\tag{\theLem}}
\newcommand{\Lhat}{\hat{L}} 
\def\:{\colon} 

\newcommand{\Zset}{\mathbb{Z}}
\newcommand{\Nset}{\mathbb{N}}
\begin{document}


\title[Homotopy fixed points for $L_{K(n)}(E_n \wedge X)$]{Homotopy fixed points for $L_{K(n)}(E_n \wedge X)$ using 
the continuous action}
\author[Daniel G. Davis]{Daniel G. Davis$\sp *$}

\date{January 12, 2005}
\maketitle

\footnotetext[1]{The author was supported by NSF grant 
DMS-9983601.}

\begin{abstract}
Let $G$ be a closed subgroup of $G_n$, the extended Morava 
stabilizer group. Let $E_n$ be the Lubin-Tate spectrum, $X$ an 
arbitrary spectrum with trivial $G$-action, and let 
$\Lhat=L_{K(n)}$. We prove that $\Lhat(E_n \wedge X)$ is a 
continuous $G$-spectrum with homotopy fixed point spectrum 
$(\Lhat(E_n \wedge X))^{hG},$ defined with respect to the 
continuous action. Also, we construct a descent spectral 
sequence whose abutment is $\pi_\ast ((\Lhat(E_n \wedge X))^{hG}).$ 
We show that the homotopy fixed points of $\Lhat(E_n \wedge X)$ 
come from the $K(n)$-localization of the homotopy fixed points of 
the spectrum $(F_n \wedge X)$.
\end{abstract}

\section{Introduction}\label{first}
Let $E_n$ be the Lubin-Tate spectrum with 
$E_{n \ast}=W\lbr u_1,...,u_{n-1} \rbr [u^{\pm 1}]$, where the 
degree of $u$ is $-2$, and
the complete power series ring over the Witt vectors 
$W =W(\mathbb{F}_{p^n})$ is in degree
zero. Let $S_n$ denote the $n$th Morava stabilizer group, the 
automorphism group of the Honda
formal group law $\Gamma_n$ of height $n$ over $\mathbb{F}_{p^n}$.
Then let $G_n=S_n\rtimes\mathrm{Gal},$ where $\mathrm{Gal}$ is
the Galois group $\mathrm{Gal}(\mathbb{F}_{p^n}/{\mathbb{F}_p})$, 
and let $G$ be
a closed subgroup of $G_n$. Morava's change of rings theorem yields 
a spectral sequence
\begin{equation}\label{ss}\zig
H^\ast_c(G_n;\pi_\ast(E_n \wedge X))\Rightarrow\pi_\ast \Lhat 
(X),\end{equation}
where the $E_2$-term is continuous cohomology, $X$ is a finite
spectrum, and $\Lhat  $ is Bousfield localization with respect
to $K(n)_\ast = \mathbb{F}_p[v_n, v_n^{-1}],$ where $K(n)$ is
Morava $K$-theory (see \cite[Prop. 7.4]{HMS}, \cite{Morava}). 
Using the $G_n$-action on $E_n$ by maps of commutative
$S$-algebras (work of Goerss and Hopkins (\cite{AndreQuillen},
\cite{Pgg/Hop0}), and Hopkins and Miller \cite{Rezk}), Devinatz and 
Hopkins \cite{DH} construct
spectra $E_n^{\mathtt{h}G}$ with strongly convergent spectral 
sequences 
\begin{equation}\label{lastlabel}\zig
H^s_c(G;\pi_t(E_n \wedge X)) \Rightarrow 
\pi_{t-s}(E_n^{\mathtt{h}G}\wedge
X).\end{equation} Also, they show that $E_n^{\mathtt{h}G_n} 
\simeq \Lhat(S^0),$ so that $E_n^{\mathtt{h}G_n} \wedge X 
\simeq \Lhat (X).$ 
\par
We recall from \cite[Rk. 1.3]{DH} the definition of the continuous 
cohomology that appears above (see also Lemma \ref{help2}). 
Let $I_n = (p, u_1, ..., u_{n-1}) 
\subset 
E_{n \ast}$. The isomorphism $\pi_t(E_n \wedge X) = 
\mathrm{lim} \, _k \, 
\pi_t(E_n \wedge X)/{I_n^k\pi_t(E_n \wedge X)}$ presents 
$\pi_t(E_n \wedge X)$ as the inverse limit of finite discrete 
$G$-modules. Then the above continuous cohomology is defined by 
\[H^s_c(G;\pi_t(E_n \wedge X)) = \mathrm{lim} \, _k \, 
H^s_c(G;\pi_t(E_n \wedge X)/{I_n^k\pi_t(E_n \wedge X)}).\]
\par
We compare the spectrum $E_n^{\mathtt{h}G}$ and spectral 
sequence (\ref{lastlabel}) with constructions for homotopy 
fixed point spectra. 
When $K$ is a discrete group and $Y$ is a $K$-spectrum of 
topological spaces, there is
a homotopy fixed point spectrum $Y^{hK}=\mathrm{Map}_K(EK_+, Y),$
where $EK$ is a free contractible $K$-space. Also, there is a 
descent spectral sequence 
\[E_2^{s,t} =
H^s(K;\pi_t(Y)) \Rightarrow \pi_{t-s}(Y^{hK}),\]
where the $E_2$-term is group cohomology \cite[\S 1.1]{Mitchell}. 
\par
Now let $K$ be a profinite group. If $S$ is a $K$-set, then $S$ 
is a {\em discrete $K$-set\/} if the action map $K \times S 
\rightarrow S$ is continuous, where $S$ is given the discrete 
topology. Then, a {\em discrete $K$-spectrum\/} $Y$ is a 
$K$-spectrum of simplicial sets, such that each simplicial set 
$Y_k$ is a simplicial discrete $K$-set (that is, for each 
$l \geq 0$, $Y_{k,l}$ is a discrete $K$-set, and all the face 
and degeneracy maps are $K$-equivariant). Then, due to work 
of Jardine and Thomason, as explained 
in Sections \ref{disdef} and \ref{hfps}, there 
is a homotopy fixed point spectrum $Y^{hK}$ defined with respect 
to the continuous action of $K$, and, in nice situations, a descent 
spectral sequence \[H^s_c(K; \pi_t(Y)) \Rightarrow 
\pi_{t-s}(Y^{hK}),\] where the $E_2$-term is the continuous 
cohomology of $K$ with coefficients in the discrete $K$-module 
$\pi_t(Y)$. 
\par
Notice that we use the notation $E_n^{\mathtt{h}G}$ for the 
construction of Devinatz and Hopkins, and $(-)^{hK}$ for homotopy 
fixed points with respect to a continuous action, although 
henceforth, when $K$ is finite and $Y$ is a $K$-spectrum of 
topological spaces, we write $Y^{h'K}$ for $\mathrm{holim} 
\, _K \, Y$, which is an equivalent definition of the homotopy 
fixed point spectrum $\mathrm{Map}_K(EK_+, Y)$. 
\par
After comparing the spectral sequence for $E_n^{\mathtt{h}G} 
\wedge X$ with the descent spectral sequence for $Y^{hK}$, $E_n 
\wedge X$ appears to be a continuous
$G_n$-spectrum with ``descent'' spectral sequences for 
``homotopy fixed point spectra'' $E_n^{\mathtt{h}G}\wedge X.$ 
Indeed, we apply \cite{DH} to show that $E_n \wedge X$ is a 
continuous
$G_n$-spectrum; that is, $E_n \wedge X$ is the homotopy limit of 
a tower of fibrant discrete
$G_n$-spectra. Using this continuous action, we define the 
homotopy
fixed point spectrum $(E_n\wedge X)^{hG}$ and construct its 
descent spectral sequence.
\par
In more detail, the $K(n)$-local spectrum $E_n$ has an
action by the profinite group $G_n$ through maps of commutative 
$S$-algebras. The spectrum $E_n^{\mathtt{h}G}$, a $K(n)$-local 
commutative
$S$-algebra, is referred to as a ``homotopy 
fixed point spectrum'' because it has the following desired 
properties: (a) for any
finite spectrum $X$, there exists spectral sequence 
(\ref{lastlabel}), 
which has the form of a descent spectral
sequence;  (b) when $G$ is finite, there is a weak 
equivalence $E_n^{\mathtt{h}G} \rightarrow E_n^{h'G},$ and 
the descent spectral sequence for $E_n^{h'G}$ is
isomorphic to spectral sequence (\ref{lastlabel}) (when $X=S^0$)
\cite[Thm. 3]{DH}; and (c) $E_n^{\mathtt{h}G}$ is an 
$N(G)/G$-spectrum, where $N(G)$ is the normalizer of $G$ in $G_n$
\cite[pg. 5]{DH}. Thus, their constructions strongly 
suggest that $G_n$ acts
on $E_n$ in a continuous sense. 
\par
However, in \cite{DH}, the action of $G_n$ on $E_n$ 
is not proven to be continuous, and $E_n^{\mathtt{h}G}$ is 
not defined with respect to a continuous $G$-action. Also, when 
$G$ is profinite, homotopy fixed points should always be the total 
right derived functor of fixed points, in some sense, and, 
in \cite{DH}, it is not shown that the ``homotopy fixed point 
spectrum'' $E_n^{\mathtt{h}G}$ can be obtained through such 
a total right derived functor.
\par
After introducing some notation, we state the main results of
this paper. Let $BP$ be the Brown-Peterson spectrum with 
$BP_\ast =
\Zset_{(p)}[v_1, v_2, ...]$, where the degree of $v_i$ 
is $2(p^i-1).$ The
ideal $(p^{i_0}, v_1^{i_1},
..., v_{n-1}^{i_{n-1}}) \subset BP_\ast$ is denoted by $I$; 
$M_I$ is the corresponding generalized Moore
spectrum $M(p^{i_0}, v_1^{i_1}, ..., v_{n-1}^{i_{n-1}}),$ a 
trivial $G_n$-spectrum. Given
an ideal $I$, $M_I$ need not exist; however, enough exist for 
our
constructions. The spectrum $M_I$ is a finite type $n$ spectrum 
with $BP_\ast(M_I) \cong BP_\ast/I.$ The set 
$\{i_0, ..., i_{n-1}\}$ of
superscripts varies so
that there is a family of ideals $\{I\}$. 
(\cite[\S 4]{Moravamodules}, 
\cite[\S4]{HS}, and
\cite[Prop. 3.7]{Mahowald/Sadofsky} provide details for our 
statements about the spectra $M_I$.) The map $r \:
BP_\ast \rightarrow E_{n \ast}$ - defined by $r(v_i) =
u_iu^{1-p^i}$, where $u_n=1$ and $u_i = 0$, when $i > n$ - makes 
$E_{n \ast}$ a $BP_\ast$-module. Then, by the Landweber exact 
functor theorem for $BP$, $\pi_\ast(E_n \wedge M_I) \cong 
E_{n \ast}/I.$ 
\par
The collection $\{I\}$ contains a descending chain of ideals 
$\{I_0 \supset I_1 \supset I_2 \supset \cdots \}$, such that 
there exists a corresponding tower of generalized Moore spectra 
$\{M_{I_0} \leftarrow M_{I_1} \leftarrow M_{I_2} \leftarrow 
\cdots\}.$ In what follows, the functors $\mathrm{lim} \, _I$ 
and
$\mathrm{holim} \, _I$ are always taken over the tower of ideals 
$\{I_i\},$ so that $\mathrm{lim} \, _I$ and $\mathrm{holim} \, _I$ 
are really $\mathrm{lim} \, _{I_i}$ and $\mathrm{holim} \, _{I_i},$ 
respectively. Also, in this paper, the homotopy limit of spectra, 
holim, is constructed levelwise in $\mathcal{S}$, the category
of simplicial sets, as defined in \cite{Bousfield/Kan} and 
\cite[5.6]{Thomason}.
\par
As in \cite[(1.4)]{DH}, let $\GG = U_0 \gneq U_1 \gneq \cdots 
\gneq 
U_i \gneq \cdots$ be a descending chain of open normal 
subgroups, such that
$\bigcap \hspace{.25pt} _{i} \, U_i = \{e\}$ and the 
canonical map
$\GG \rightarrow \mathrm{lim} \, _i \, \GG/{U_i}$ is a 
homeomorphism. We define
\[F_n = 
\mathrm{colim} \, _i \, E_n^{\mathtt{h}U_i}.\] 
Then the key to getting our work started is knowing that
\[E_n \wedge M_I \simeq F_n \wedge M_I,\] and thus, $E_n \wedge
M_I$ has the homotopy type of the discrete
$G_n$-spectrum $F_n \wedge M_I.$ This result (Corollary 
\ref{Rezk}) 
is not difficult, thanks to the work of Devinatz and Hopkins.
\par
Given a tower $\{Z_I\}$ of discrete $G_n$-spectra, there is a 
tower 
$\{(Z_I)_f\},$ with $G_n$-equivariant maps $Z_I \rightarrow 
(Z_I)_f$ 
that are weak equivalences, and $(Z_I)_f$ is a fibrant discrete 
$G_n$-spectrum (see Def. \ref{fibrep}). For the remainder 
of this section, $X$ is any spectrum with trivial $G_n$-action. 
We use $\cong$ to denote an isomorphism in the stable homotopy 
category.
\begin{Thm}\label{theorem}\label{m1}
As the homotopy limit of a tower of fibrant discrete $G_n$-spectra, 
$E_n \cong \mathrm{holim} \, _I \, (F_n \wedge
 M_I)_f$ is a continuous $G_n$-spectrum. Also, for any 
spectrum $X$, $\Lhat  (E_n \wedge X) \cong
\mathrm{holim} \, _I \, (F_n \wedge
 M_I \wedge X)_f$ is a continuous $G_n$-spectrum.
\end{Thm}  
\par
Using the hypercohomology spectra of Thomason \cite{Thomason}, and 
the theory of pre- sheaves of spectra on a site and globally 
fibrant 
models, developed by Jardine (e.g. \cite{Jardinejpaa}, 
\cite{Jardinecanada}, 
\cite{Jardine}, \cite{Jardinesummary}), we present the 
theory of homotopy fixed points 
for discrete $G$-spectra by considering the site of finite discrete 
$G$-sets. Using the fact that $G_n$ has finite virtual cohomological 
dimension, Thomason's hypercohomology spectrum gives a concrete 
model for the homotopy fixed points that makes building the descent 
spectral sequence easy. We apply this theory to define homotopy 
fixed points for towers 
of discrete $G$-spectra.
\par
We point out that much of the theory described above (in 
Sections 3, 5, 7, and 8, through Remark \ref{l-adic}) is 
already known, in some form, especially in the work of Jardine 
mentioned above, in the excellent article \cite{Mitchell}, 
by Mitchell (see also the opening remark of 
\cite[\S 5]{hypercohomology}), and in Goerss's paper \cite{hGal}. 
However, since the above theory has not been explained in detail 
before, using the language of homotopy fixed points for discrete 
$G$-spectra, we give a presentation of it.
\par

After defining homotopy fixed points for towers of discrete 
$G$-spectra, we show that these homotopy fixed points are the 
total right derived functor of fixed points in the appropriate 
sense, and we construct the associated descent spectral sequence. 
This enables us to define 
the homotopy fixed point spectrum $(\Lhat  (E_n \wedge X))^{hG}$, 
using the continuous $G_n$-action, and construct its descent 
spectral sequence. More specifically, we have the following 
results.
\begin{Def}
Given a profinite group $G$, let $\mathcal{O}_G$ be the {\em 
orbit category\/} of $G$. The objects of $\mathcal{O}_G$ are the 
continuous 
left $G$-spaces $G/K$, for all $K$ closed in $G$, and the 
morphisms are the 
continuous $G$-equivariant maps. Note that each object in 
$\mathcal{O}_G$ is a profinite space.
\end{Def}
\par
We use ${Sp}$ to denote the model category 
$\mathrm{(spectra)}^\mathrm{stable}$ of Bousfield-Friedlander 
spectra.
\begin{Thm}\label{m2}
There is a functor \[P \: (\mathcal{O}_{G_n})^\mathrm{op} 
\rightarrow Sp, \ \ \ P(G_n/G) = E_n^{hG},\] where $G$ is any 
closed subgroup of $G_n$.
\end{Thm}
\par
We also show that the $G$-homotopy fixed points of 
$\Lhat(E_n \wedge X)$ can be obtained by taking the 
$K(n)$-localization of the $G$-homotopy fixed points of the 
discrete $G$-spectrum $(F_n \wedge X)$. This result shows that 
the spectrum $F_n$ is an interesting spectrum that is worth 
further study.  
\begin{Thm}\label{m3}
For any closed subgroup $G$ and any spectrum $X$, there is an 
isomorphism \[(\Lhat(E_n \wedge X))^{hG} \cong \Lhat 
((F_n \wedge X )^{hG})\] in the stable homotopy category. 
In particular, $E_n^{hG} \cong \Lhat(F_n^{hG})$.
\end{Thm}
\begin{Thm}\label{m4}
Let $G$ be a closed subgroup of $G_n$ and let $X$ be any spectrum. 
Then there is a conditionally
convergent descent spectral sequence 
\begin{equation}\label{nice2}\zig E_2^{s,t} \Rightarrow 
\pi_{t-s}((\Lhat  (E_n \wedge X))^{hG}).\end{equation} 
If the tower of abelian groups $\{\pi_t(E_n \wedge 
M_I \wedge X)\}_I$ 
satisfies the Mittag-Leffler condition, for each $t \in 
\Zset$, 
then \[E_2^{s,t} \cong H^s_\mathrm{cts}(G;\pi_t(\Lhat (E_n \wedge
X))),\] the cohomology of continuous cochains. If $X$ is a finite 
spectrum, then (\ref{nice2}) has the form
\begin{equation}\label{nice}\zig H^s_c(G;
\pi_t(E_n \wedge X)) \Rightarrow \pi_{t-s}((E_n \wedge
X)^{hG}),\end{equation} where the $E_2$-term is the continuous 
cohomology of $\mathrm{(}\ref{lastlabel}\mathrm{)}$.
\end{Thm}
\par
Also, Theorem \ref{final} shows that, when $X$ is finite, 
$(E_n \wedge X)^{hG} \cong E_n^{hG} \wedge X$, so that descent 
spectral sequence (\ref{nice}) has the same form as spectral 
sequence 
(\ref{lastlabel}). It is natural to wonder if these two spectral 
sequences are isomorphic to each other. Also, the spectra 
$E_n^{\mathtt{h}G}$ and $E_n^{hG}$ should be the same. We plan 
to say more about the relationship between 
$E_n^{\mathtt{h}G}$ and $E_n^{hG}$ and their associated spectral 
sequences in future work. 
\par
It is important to note that the model given here for $E_n$ as
a continuous $\GG$-spectrum is not completely satisfactory, since 
the continuous action is by morphisms
that are just maps of spectra. Because there are models for $E_n$ 
where the $\GG$-action is by
$A_\infty$- and $E_\infty$-maps of ring
spectra, we would like to know that such structured actions are
actually continuous.
\vspace{.1in}
\par
We outline the contents of this paper. In Section 2, after 
establishing some notation and terminology, we provide some 
background material and recall useful facts. In Section 3, we 
study the model category of discrete $G$-spectra. In Section 4, 
we study towers of discrete $G$-spectra and we give a 
definition 
of continuous $G$-spectrum. In Section 5, we define homotopy 
fixed 
points for discrete $G$-spectra and state some basic facts about 
this 
concept. Section 6 shows that $E_n$ is a continuous 
$G_n$-spectrum, 
proving the first half of Theorem \ref{m1}. Section 7 constructs 
two useful models of the $G$-homotopy fixed point spectrum, 
when $G$ has finite virtual cohomological dimension. Section 8 
defines homotopy fixed points for towers of discrete $G$-spectra, 
builds a descent spectral sequence in this setting, and 
shows that 
these homotopy fixed points are a total right derived functor, 
in 
the appropriate sense. Section 9 
completes the proof of Theorem 
\ref{m1}, studies $(\Lhat(E_n \wedge X))^{hG}$, 
and proves Theorems \ref{m2} and \ref{m3}. Section 10 considers 
the 
descent spectral sequence for $(\Lhat(E_n \wedge X))^{hG}$ and 
proves 
Theorem \ref{m4}.

\vspace{.1in}
\textbf{Acknowledgements.} This paper is a development of part 
of my thesis. I am 
very grateful to my thesis advisor, Paul Goerss, for many helpful 
conversations and useful suggestions regarding this paper. 
Also, I thank Ethan Devinatz for 
very helpful answers to my questions about his work \cite{DH} 
with Mike Hopkins. I am grateful to Halvard Fausk, Christian 
Haesemeyer, Rick 
Jardine, and Charles Rezk for useful conversations.

\section{Notation, Terminology, and Preliminaries}\label{prelims} 
\par
We begin by establishing some notation and terminology that 
will be 
used throughout the paper. $\mathbf{Ab}$ is the category of 
abelian groups. Outside of $\mathbf{Ab}$, all groups are assumed 
to be 
profinite,
unless stated otherwise. For a group $G$, we write $G \cong
\mathrm{lim} \, _{N} \, G/{N}$, the inverse limit over the open 
normal subgroups. The notation $H < _c G$ means that $H$ is
a closed subgroup of $G$. We use $G$ to denote arbitrary 
profinite groups and,
specifically, closed subgroups of $G_n$.
\par
Let $\mathcal{C}$ be a category. A tower $\{C_i\}$ of objects in 
$\mathcal{C}$ is a diagram in $\mathcal{C}$ of the form $\cdots 
\rightarrow C_i \rightarrow C_{i-1} \rightarrow \cdots \rightarrow 
C_1 \rightarrow C_0.$ We always use Bousfield-Friedlander spectra 
\cite{BF}, except 
when 
another category of spectra is specified. If $\mathcal{C}$ is a 
model
category, then $\mathrm{Ho}(\mathcal{C})$ is its homotopy
category. The phrase ``stable category'' always refers to
$\mathrm{Ho}(Sp)$.
\par
In $\mathcal{S}$, the category of simplicial sets, $S^n = 
\Delta^n/{\partial \Delta^n}$ is the $n$-sphere. Given a spectrum 
$X$, $X^{(0)}=S^0$, and for $j\geq 1$, $X^{(j)} = X \wedge X \wedge 
\cdots 
\wedge X,$ with
$j$ factors. $L_n$ denotes 
Bousfield localization with respect to $E(n)_\ast 
\negthinspace = \negthinspace
\Zset_{(p)}[v_1, ..., v_{n-1}][v_n,v_n^{-1}]$. 
\begin{Def}
\cite[Def. 1.3.1]{Hoveybook} Let $\mathcal{C}$ and $\mathcal{D}$ 
be model categories. The
functor $F \: \mathcal{C} \rightarrow \mathcal{D}$ is a {\em 
left
Quillen functor\/} if $F$ is a left adjoint that preserves
cofibrations and trivial cofibrations. The functor $P \:
\mathcal{D} \rightarrow \mathcal{C}$ is a {\em right Quillen 
functor\/}
if $P$ is a right adjoint that preserves fibrations and trivial
fibrations. Also, if $F$ and $P$ are an adjoint pair and left and 
right Quillen functors, respectively, then $(F,P)$ is a {\em 
Quillen pair\/} for the model categories
$(\mathcal{C}, \mathcal{D})$. 
\end{Def}
Recall \cite[Lemma 1.3.10]{Hoveybook} that a Quillen pair
$(F,P)$ yields total left and right derived functors
$\mathbf{L}F$ and $\mathbf{R}P$, respectively, which give an 
adjunction between the homotopy categories 
$\mathrm{Ho}(\mathcal{C})$ and
$\mathrm{Ho}(\mathcal{D})$. 
\par 
We use $\mathrm{Map}_c(G, A) = \Gamma_G(A)$ to 
denote the set of continuous maps from $G$ to the topological 
space $A$, where $A$ is often a set, equipped with the discrete 
topology, or a discrete abelian group. Instead of $\Gamma_G (A)$, 
sometimes we write just $\Gamma (A)$, when the $G$ is understood 
from context. Let $(\Gamma_G)^k (A)$ denote $(\Gamma_G \, \Gamma_G 
\cdots \Gamma_G)(A)$, the application of $\Gamma_G$ to $A$, 
iteratively, $k+1$ times, where $k \geq 0$. Let $G^k$ be the 
$k$-fold product of $G$ and let $G^0=\ast$. Then, if $A$ is a 
discrete set or a 
discrete abelian group, there is a $G$-equivariant isomorphism 
$(\Gamma_G)^k(A) \cong \mathrm{Map}_c(G^{k+1},A)$ of discrete 
$G$-sets 
(modules), where $\mathrm{Map}_c(G^{k+1},M)$ has $G$-action 
defined 
by \[(g' \cdot
f)(g_1, ..., g_{k+1}) = f(g_1g', g_2, g_3, ..., g_{k+1}).\] 
Also, we often 
write $\Gamma_G^k(A)$, or $\Gamma^k A$, for 
$\mathrm{Map}_c(G^k,A).$ 
\par
Let $A$ be a discrete abelian group. Then 
$\mathrm{Map}_c^\ell(G_n^k,A)$ is
the discrete $G_n$-module of continuous maps $G_n^k \rightarrow A$ 
with action defined by \[(g' \cdot
f)(g_1, ..., g_k) = f((g')^{-1}g_1, g_2, g_3, ..., g_k).\] 
It is helpful to note that there is a $G_n$-equivariant 
isomorphism 
of discrete $G_n$-modules \[p \: \mathrm{Map}_c^\ell(G_n^k,A) 
\rightarrow
\mathrm{Map}_c(G_n^k,A), \ \ p(f)(g_1, g_2, ..., g_k) =
f(g_1^{-1},g_2, ..., g_k).\]
$\mathrm{Map}_c^\ell(G_n^k,A)$ is also defined when $A$ is an 
inverse
limit of discrete abelian groups.
\par
By a {\em topological $G$-module}, we mean an abelian Hausdorff 
topological group that is a $G$-module, with a continuous 
$G$-action. 
Note that if $M = \mathrm{lim} \, _i \, M_i$ is the inverse limit 
of a tower $\{M_i\}$ of discrete $G$-modules, then $M$ is a 
topological $G$-module.
\par
For the remainder of this section, we recall some frequently used 
facts and discuss
background material, to help
get our work started.
\par
In \cite{DH}, Devinatz and Hopkins, using work by
Goerss and Hopkins (\cite{AndreQuillen},
\cite{Pgg/Hop0}), and Hopkins and Miller \cite{Rezk}, show that 
the action of $\GG$ on $E_n$ is by maps
of commutative $S$-algebras. Previously, Hopkins and Miller
had shown that $G_n$ acts on $E_n$ by maps of $A_\infty$-ring
spectra. However, the continuous action presented here is not 
structured. As already mentioned, the starting point for the
continuous action is the spectrum $F_n \wedge M_I$, which is not 
known to be an $A_\infty$-ring object in the category of discrete 
$G_n$-spectra. Thus, we work in the unstructured category $Sp$ of 
Bousfield-Friedlander spectra of
simplicial sets, and the continuous action is simply by maps of 
spectra.
\par
As mentioned above, \cite{DH} is written using $E_\infty$, the 
category of commutative $S$-algebras, and $\mathcal{M}_S$, the 
category of $S$-modules (see \cite{EKMM}). However, 
\cite[\S 4.2]{HSS}, 
\cite[\S 14, \S 19]{mmss}, and \cite[pp. 529-530]{Schwede} 
show that $\mathcal{M}_S$ and $Sp$ are Quillen equivalent model 
categories \cite[\S 1.3.3]{Hoveybook}. Thus, we can import the 
results of Devinatz and Hopkins from $\mathcal{M}_S$ into $Sp$. 
For example, \cite[Thm. 1]{DH} implies the following result, 
where $R^+_{\GG}$ is the category whose
objects are finite discrete left $\GG$-sets and $G_n$ itself 
(a continuous
profinite left $\GG$-space), and whose morphisms are continuous
$\GG$-equivariant maps. 
\begin{Thm}[Devinatz, Hopkins]
There is a presheaf of spectra \[F \colon 
(R^+_{\GG})^{\mathrm{op}} 
\rightarrow Sp,\] such that (a) for each $S \in R^+_{\GG}$,
$F(S)$ is $K(n)$-local; (b) $F(\GG) = E_n$; (c) for $U$ an 
open subgroup of 
$\GG$, 
$E_n^{\mathtt{h}U} := F(\GG/U)$;
and (d) $F(\ast) = E_n^{\mathtt{h}\GG} \simeq \Lhat  S^0.$ 
\end{Thm}
\par
Now we define a spectrum that is essential to our
constructions.
\begin{Def}\label{center} Let $F_n = \mathrm{colim} \, _i \,
E_n^{\mathtt{h}U_i},$ where the direct limit is in $Sp$. Because 
$\mathrm{Hom}_{G_n}(G_n/{U_i}, G_n/{U_i}) \cong G_n/{U_i},$ $F$ 
makes $E_n^{\mathtt{h}U_i}$ a $G_n/{U_i}$-spectrum. Thus, $F_n$ 
is a $G_n$-spectrum, and the canonical map $\eta \: F_n 
\rightarrow 
E_n$ is $G_n$-equivariant.
\end{Def}
\par
Note that $F_n$ is the stalk of the presheaf of spectra
${F}|_{(G\negthinspace - \negthinspace 
\mathbf{Sets}_{df})^\mathrm{op}}$, at the unique point of the 
Grothendieck topos (see \S \ref{spt}). The following useful fact 
is stated in \cite[pg. 9]{DH} (see 
also \cite[Lemma 14]{Str}).
\begin{Thm}\label{powers}
For $j \geq 0,$ let $X$ be a finite spectrum and regard 
$\Lhat(E_n^{(j+1)} \wedge X)$ as a $G_n$-spectrum, where $G_n$ 
acts only on the leftmost factor of the smash product. Then 
there is a $G_n$-equivariant isomorphism 
\[\pi_\ast(\Lhat(E_n^{(j+1)} \wedge X)) \cong 
\mathrm{Map}^\ell_c(G_n^j,\pi_\ast(E_n \wedge X)).\] 
\end{Thm}
\par
We review some frequently used facts about the functor $L_n$ 
and homotopy limits of spectra. First, $L_n$ 
is smashing, e.g. $L_n X \simeq X \wedge L_n S^0,$ for any 
spectrum $X$, and $E(n)$-localization commutes with homotopy 
direct limits
\cite[Thms. 7.5.6, 8.2.2]{Ravenelorange}. Note that this implies 
that $F_n$ is $E(n)$-local.
\begin{Def}
If $\ \cdots \rightarrow X_i \rightarrow X_{i-1} \rightarrow 
\cdots \rightarrow X_1 \rightarrow X_0$  is a tower of spectra 
such that each $X_i$ is fibrant in $Sp$, then $\{X_i\}$ is a 
{\em tower of fibrant spectra}. 
\end{Def}
\par
If $\{X_i\}$ is a tower of fibrant spectra, then there is a short 
exact sequence \[0 \rightarrow \textstyle{\mathrm{lim}}^1 \, _i 
\, \pi_{m+1}(X_i) \rightarrow \pi_m(\mathrm{holim} \, _i \, X_i) 
\rightarrow \mathrm{lim} \, _i \, \pi_m(X_i) \rightarrow 0.\] 
Also, if each map in the tower is a fibration, the map 
$\mathrm{lim} \, _i \, X_i \rightarrow \mathrm{holim} \, _i \, 
X_i$ is a weak equivalence. If $J$ is a small category and the 
functor $P \: J \rightarrow Sp$ is a diagram of spectra, such 
that $P_j$ is fibrant for each $j \in J$, then $\mathrm{holim} 
\, _j \, P_j$ is a fibrant spectrum.
\begin{Def}
There is a functor $(-)_\mathtt{f} \: Sp \rightarrow Sp,$ such 
that, given $Y$ in $Sp$, $Y_\mathtt{f}$ is a fibrant spectrum, 
and there is a natural transformation $\mathrm{id}_{Sp} 
\rightarrow (-)_\mathtt{f}$, such that, for any $Y$, the map $Y 
\rightarrow Y_\mathtt{f}$ is a trivial cofibration. 
For example, if $Y$ is a $G$-spectrum, 
then $Y_\mathtt{f}$ is also a $G$-spectrum, and the map $Y 
\rightarrow Y_\mathtt{f}$ is $G$-equivariant.
\end{Def}
\par
The following statement says that smashing with a finite 
spectrum commutes with homotopy limits.
\begin{Lem}[{\cite[pg. 96]{topos}}]\label{commute}
Let $J$ be a small category, $\{Z_j\}$ a $J$-shaped diagram 
of fibrant spectra, and let $Y$ be a finite spectrum. Then 
the composition \[(\mathrm{holim} \, _j \,
Z_j) \wedge Y \rightarrow \mathrm{holim} \, _j \, (Z_j \wedge Y) 
\rightarrow \mathrm{holim} \, _j \, (Z_j \wedge
Y)_\mathtt{f}\] is a weak equivalence.
\end{Lem}
\par
We recall the result that is used to build towers of discrete 
$G$-spectra.
\begin{Thm}[{\cite[\S 2]{HoveyCech}, 
\cite[Remark 3.6]{DevFields}}]\label{veryneat}
If $X$ is an $E(n)$-local spectrum, then, in the stable category, 
there is an isomorphism \[\Lhat  X \cong
\mathrm{holim} \, _I \, (X \wedge M_I)_\mathtt{f}.\]
\end{Thm}
\begin{Lem}[{\cite[Lemma 7.2]{HS}}]
If $X$ is any spectrum, and $Y$ is a finite spectrum
of type $n$, then $\Lhat  (X
\wedge Y) \simeq \Lhat  (X)
\wedge Y \simeq L_n(X) \wedge Y.$
\end{Lem}
\par
We recall some useful facts about compact $p$-adic analytic
groups. Since $S_n$ is compact
$p$-adic analytic, and $\GG$
is an extension of ${S}_n$ by
the Galois group, $\GG$ is a compact $p$-adic analytic group 
\cite[Cor. of Thm. 2]{Serrep}. Any closed subgroup of a compact 
$p$-adic analytic group is also compact $p$-adic
analytic \cite[Thm. 9.6]{Dixon}. Also, since the subgroup in 
$S_n$ of strict automorphisms is finitely generated and pro-$p$, 
\cite[pp. 76, 124]{Ribes} implies that all subgroups in $G_n$ of 
finite index are open. 
\par
Let the profinite group $G$ be a compact $p$-adic 
analytic group. Then $G$ contains an open subgroup $H$, 
such that $H$ is a pro-$p$ group with finite cohomological 
$p$-dimension; that is, $\mathrm{cd}_p(H)=m$, for some 
non-negative integer $m$ (see \cite[2.4.9]{Lazard} or the 
exposition in \cite{Symonds}). Since $H$ is pro-$p$, 
$\mathrm{cd}_q(H)=0$, whenever $q$ is a prime different 
from $p$ \cite[Prop. 11.1.4]{Wilson}. Also, if $M$ is a 
discrete $H$-module, then, for $s \geq 1$, $H^s_c(H;M)$ is a 
torsion abelian group \cite[Cor. 6.7.4]{Ribes}. These facts 
imply that, for any discrete $H$-module $M$, $H^s_c(H;M)=0$, 
whenever $s>m+1$. We express this conclusion by saying that 
$G$ has finite virtual cohomological dimension and we write 
$\mathrm{vcd}(G) \leq m$. Also, if $K$ is a closed subgroup of 
$G$, $H \cap K$ is an open pro-$p$ subgroup of $K$ with 
$\mathrm{cd}_p(H \cap K) \leq m$, so that 
$\mathrm{vcd}(K) \leq m$, and thus, $m$ is a uniform 
bound independent of $K$.
\par
Now we state various results related to towers of abelian groups 
and 
continuous cohomology. The lemma below follows from the fact that 
a tower 
of abelian groups satisfies the Mittag-Leffler condition if and 
only if the tower is pro-isomorphic to a tower of epimorphisms 
\cite[(1.14)]{Jannsen}. 
\begin{Lem}\label{mittag} 
Let $F \: \mathbf{Ab} \rightarrow \mathbf{Ab}$ be an exact 
additive 
functor. If \ $\{A_i\}_{i \geq 0}$ is a tower of abelian groups 
that satisfies the Mittag-Leffler condition, then so does the 
tower 
$\{F(A_i)\}.$
\end{Lem}
\begin{Rk}\label{mittagrk}
Let $G$ be a profinite group. Then the functor 
\[\mathrm{Map}_c(G,-) \: \mathbf{Ab} \rightarrow \mathbf{Ab}, 
\ \ \ A \mapsto \mathrm{Map}_c(G,A),\] is defined by giving 
$A$ the discrete topology. The isomorphism 
$\mathrm{Map}_c(G,A) \cong \mathrm{colim} \, _N \prod 
\, _{G/N} \, A$ shows that $\mathrm{Map}_c(G,-)$ is an 
exact additive functor. Later, we will use Lemma \ref{mittag} 
with this functor.
\end{Rk}
\par
The next lemma is a consequence of the fact that limits 
in 
$\mathbf{Ab}$ and in topological spaces are created in 
$\mathbf{Sets}$. 
\begin{Lem}
Let $M = \mathrm{lim} \, _\alpha \, M_\alpha$ be an inverse limit
of discrete abelian groups, so that $M$ is an abelian 
topological group. Let $H$ be any profinite group. 
Then $\mathrm{Map}_c(H,M) \rightarrow \textstyle{\mathrm{lim}} 
\, _\alpha \,
\mathrm{Map}_c(H, M_\alpha)$ is an isomorphism of abelian groups.
\end{Lem}
\begin{Lem}\label{finite}
If $X$ is a finite spectrum, $G<_c G_n$, and $t$ any integer, 
then the abelian group $\pi_t(E_n \wedge M_I \wedge X)$ is finite.
\end{Lem}
\begin{proof}
The starting point is the fact that $\pi_t(E_n \wedge M_I) 
\cong \pi_t(E_n)/I$ is finite. In the stable homotopy category of
CW-spectra, since $X$ is a finite spectrum, there exists 
some $m$ such that $X_l = \Sigma^{l-m} 
X_m$
whenever $l \geq m,$ and $X_m$ is a finite complex. Since
$\pi_t(E_n \wedge M_I \wedge X) \cong 
\pi_{t-m}(E_n \wedge M_I \wedge X_m)$ and $X_m$ 
can be built out of a finite number of
cofiber sequences, the result follows.
\end{proof}
\begin{Cor}[{\cite[pg. 116]{HMS}}]\label{finite2}
If $X$ is a finite spectrum, then 
$\pi_t(E_n \wedge
X) \cong \mathrm{lim} \, _I \, \pi_t(E_n \wedge M_I \wedge 
X).$
\end{Cor}
\par
We recall the definition of a second version of continuous 
cohomology. If $M$ is a topological $G$-module, then 
$H^s_\mathrm{cts}(G;M)$ is the $s$th cohomology group of the 
cochain 
complex $M \rightarrow \mathrm{Map}_c(G,M) \rightarrow 
\mathrm{Map}_c(G^2,M) \rightarrow \cdots$ of continuous cochains 
for a profinite group $G$ 
with
coefficients in $M$ (see \cite[\S2]{Tate}). If $M$ is a 
discrete $G$-module, this is the
usual continuous cohomology $H^s_c(G;M)$. There is the following 
useful relationship between 
these cohomology theories. 
\begin{Thm}[{\cite[(2.1), Thm. 2.2]{Jannsen}}]\label{ses}
Let $\{M_n\}_{n \geq 0}$ be a tower of discrete
$G$-modules satisfying the Mittag-Leffler condition and let 
$M =
\mathrm{lim} \, _n \, M_n$ as a topological $G$-module. Then, 
for
each $s \geq 0$, there is a short exact sequence \[0 \rightarrow
\textstyle{\mathrm{lim}}^1 \, _n \, H^{s-1}_c(G; M_n) 
\rightarrow H^{s}_\mathrm{cts}(G;
M) \rightarrow \mathrm{lim} \, _n \, H^{s}_c(G;
M_n) \rightarrow 0,\] where $H^{-1}_c(G;-)=0.$
\end{Thm}
\par
The next result is from \cite[Rk. 1.3]{DH} and is due to the 
fact that, when $X$ is a finite spectrum, $\pi_t(E_n \wedge X) 
\cong \mathrm{lim} \, _k \, \pi_t(E_n \wedge X)/{I_n^k\pi_t(E_n 
\wedge X)}$ is a profinite continuous $\Zset_p \lbr G 
\rbr$-module. 
\begin{Lem}\label{help2}
If $G$ is closed in $G_n$ and $X$ is a finite spectrum, then, for 
$s \geq
0,$
\begin{align*}
H^s_\mathrm{cts}(G;\pi_t(E_n \wedge X)) & \cong 
\mathrm{lim} \, _k \,
H^s_c(G; \pi_t(E_n \wedge X)/{I_n^k\pi_t(E_n \wedge X)}) \\ & = 
H^s_c(G; \pi_t(E_n \wedge X)).\end{align*}
\end{Lem}
\par
We make a few remarks about the functorial smash product in
$Sp$, defined in \cite[Chps. 1, 2]{Jardine}. 
\begin{Def}
Given spectra $X$ and $Y$, their smash product $X \wedge Y$ is
given by $(X \wedge Y)_{2k} = X_k \wedge Y_k$ and 
$(X \wedge Y)_{2k+1} = X_k \wedge Y_{k+1},$ where, for example, 
$X_k \in \mathcal{S}_\ast$, the category of pointed simplicial 
sets.
\end{Def}
\par
Since $K \wedge (-) \: \mathcal{S}_\ast \rightarrow
\mathcal{S}_\ast$ is a left adjoint, for any $K$ in
$\mathcal{S}_\ast$, smashing with any spectrum in either 
variable
commutes with colimits in ${Sp}.$  
\section{The model category of discrete $G$-spectra}\label{spt}
\par
A pointed simplicial discrete $G$-set is a pointed simplicial 
set that is a simplicial discrete $G$-set, such that the 
$G$-action fixes the basepoint. 
\begin{Def} A {\em discrete $G$-spectrum\/} $X$ is a spectrum 
of pointed simplicial sets $X_k$, for $k \geq 0$, such that each 
simplicial set $X_k$ is a pointed simplicial
discrete $G$-set, and each bonding map $S^1 \wedge X_k 
\rightarrow
X_{k+1}$ is $G$-equivariant ($S^1$ has trivial $G$-action). 
Let $Sp_G$ denote the
category of discrete $G$-spectra, where the morphisms are 
$G$-equivariant maps of spectra. 
\end{Def}
\par
As with discrete $G$-sets, if $X \in Sp_G$, there is a 
$G$-equivariant isomorphism $X \cong \mathrm{colim} \, _N \, 
X^N$. Also, a discrete
$G$-spectrum $X$ is a continuous $G$-spectrum since, for all 
$k, l \geq 0$,  the set $X_{k,l}$ is a continuous $G$-space 
with the
discrete topology, and all the face and degeneracy maps are
(trivially) continuous. 
\begin{Def}
As in \cite[\S 6.2]{Jardine}, let $G\negthinspace - 
\negthinspace 
\mathbf{Sets}_{df}$ be the canonical site of finite
discrete $G$-sets. The pretopology of $G\negthinspace - 
\negthinspace \mathbf{Sets}_{df}$ is given by covering families 
of the form $\{f_\alpha \: S_\alpha \rightarrow S \}$, a finite 
set of $G$-equivariant functions in $G\negthinspace - 
\negthinspace \mathbf{Sets}_{df}$ for a fixed $S \in G
\negthinspace - \negthinspace \mathbf{Sets}_{df}$, such that 
$\coprod \, _\alpha \, S_\alpha \rightarrow S$ is a surjection. 
\end{Def}
\par
We use $\mathbf{Shv}$ to denote the category of sheaves of sets 
on the site $G\negthinspace - \negthinspace 
\mathbf{Sets}_{df}$. Also, $T_G$ signifies the category of 
discrete $G$-sets, and, 
as in \cite{hGal}, $S_G$ is the category of simplicial 
objects in $T_G$. The Grothendieck topos $\mathbf{Shv}$ has a 
unique point 
$u \: \mathbf{Sets} \rightarrow
\mathbf{Shv}.$ The left adjoint of the topos
morphism $u$ is \[u^\ast \: \mathbf{Shv}
\rightarrow \mathbf{Sets}, \ \ \ \mathcal{F} \mapsto
\mathrm{colim} \, _N \, \mathcal{F}(G/{N}),\] with
right adjoint \[u_\ast \: \mathbf{Sets} \rightarrow \mathbf{Shv}, 
\ \ \ X \mapsto \mathrm{Hom}_G(-,
\mathrm{Map}_c(G,X))\] \cite[Rk. 6.25]{Jardine}. The $G$-action
on the discrete $G$-set $\mathrm{Map}_c(G,X)$ is defined by $(g 
\cdot f)(g') = f(g'g)$,
for $g, g'$ in $G$, and $f$ a continuous map $G \rightarrow
X$, where $X$ is given the discrete topology. 
\par
The functor $\mathrm{Map}_c(G,-) \: \mathbf{Sets} \rightarrow T_G$ 
prolongs to the functor \[\mathrm{Map}_c(G,-) \: Sp \rightarrow 
Sp_G.\] Thus, if $X$ is a spectrum, then $\mathrm{Map}_c(G,X) 
\cong \mathrm{colim} \, _N \, \prod \, _{G/N} \, X$ is the discrete 
$G$-spectrum with $(\mathrm{Map}_c(G,X))_k = 
\mathrm{Map}_c(G,X_k),$ 
where $\mathrm{Map}_c(G,X_k)$ is a pointed simplicial set, with 
$l$-simplices $\mathrm{Map}_c(G,X_{k,l})$ and basepoint $G 
\rightarrow \ast,$ where $X_{k,l}$ is regarded as a discrete set. 
The $G$-action on $\mathrm{Map}_c(G,X)$ is defined 
by the $G$-action on the sets $\mathrm{Map}_c(G,X_{k,l})$. 
\par
It is not hard to see that $\mathrm{Map}_c(G,-)$ is right adjoint 
to the forgetful functor $U \: Sp_G \rightarrow Sp$. Note that if 
$X$ is a discrete $G$-spectrum, then there is a contravariant 
functor (presheaf) $\mathrm{Hom}_G(-,X) \: (G\negthinspace - 
\negthinspace \mathbf{Sets}_{df})^\mathrm{op} \rightarrow Sp,$ 
where, for any $S \in G\negthinspace - \negthinspace 
\mathbf{Sets}_{df},$ $\mathrm{Hom}_G(S,X)$ is the spectrum with 
$(\mathrm{Hom}_G(S,X))_{k,l} = \mathrm{Hom}_G(S,X_{k,l}),$ a 
pointed set with basepoint $S \rightarrow \ast.$
\par
Let $\mathbf{ShvSpt}$ be the category of sheaves of spectra on 
the site
$G\negthinspace - \negthinspace \mathbf{Sets}_{df}$. A sheaf of 
spectra
$\mathcal{F}$ is a presheaf $\mathcal{F} \:
(G\negthinspace - \negthinspace \mathbf{Sets}_{df})^\mathrm{op} 
\rightarrow
Sp,$ such that, for any $S \in G\negthinspace - \negthinspace 
\mathbf{Sets}_{df}$ and any covering family $\{f_\alpha \: 
S_\alpha \rightarrow S \}$, the usual diagram (of spectra) is 
an equalizer. Equivalently, a sheaf of
spectra $\mathcal{F}$ consists of pointed simplicial sheaves
$\mathcal{F}^n$, together with pointed maps of simplicial
presheaves $\sigma \: S^1 \wedge
\mathcal{F}^n \rightarrow \mathcal{F}^{n+1},$ for $n \geq 0,$
where $S^1$ is the constant simplicial presheaf.
A morphism between sheaves of spectra is a natural transformation
between the underlying presheaves.
\par 
The category $\mathbf{PreSpt}$ of
presheaves of spectra on the site $G\negthinspace - \negthinspace 
\mathbf{Sets}_{df}$ has the following ``stable'' model category 
structure
(\cite{Jardinecanada}, \cite[\S 2.3]{Jardine}). A map $h \: 
\mathcal{F}
\rightarrow \mathcal{G}$ of presheaves of spectra is a weak 
equivalence
if and only if
the associated map of stalks $\colim \, _N \, \mathcal{F}(G/{N})
\rightarrow \colim \, _N \, \mathcal{G}(G/{N})$ is a weak
equivalence of spectra. Recall that a map $k$ of simplicial
presheaves is a cofibration if, for each $S \in G\negthinspace - 
\negthinspace \mathbf{Sets}_{df},$ $k(S)$
is a monomorphism of simplicial
sets. Then $h$ is a cofibration
of presheaves of spectra if the following two conditions hold:
\begin{enumerate}
\item 
the map $h^0 \: \mathcal{F}^0 \rightarrow \mathcal{G}^0$ is a 
cofibration of
simplicial presheaves; and
\item
for each $n \geq 0$, the canonical map $(S^1 \wedge 
\mathcal{G}^n) \cup_{S^1 \wedge \mathcal{F}^n}\mathcal{F}^{n+1}
\rightarrow \mathcal{G}^{n+1}$ is a cofibration of
simplicial presheaves.\end{enumerate} Fibrations are those 
maps with the right
lifting property with respect to trivial cofibrations. 
\begin{Def}
In the stable model category structure, fibrant presheaves are 
often referred to as {\em globally fibrant}, and if
$\mathcal{F} \rightarrow
\mathcal{G}$ is a weak equivalence of presheaves, with 
$\mathcal{G}$ globally fibrant, then $\mathcal{G}$ is a {\em 
globally fibrant model\/} for
$\mathcal{F}$. We often use $\mathbf{G}\mathcal{F}$ to denote 
such a globally fibrant model.
\end{Def}
\par
We recall the following fact, which is especially useful when 
$S =\ast$.
\begin{Lem}\label{Kan}
Let $S \in G\negthinspace - \negthinspace \mathbf{Sets}_{df}$. 
The $S$-sections functor
$\mathbf{PreSpt} \rightarrow
Sp$, defined by $\mathcal{F} \mapsto \mathcal{F}(S)$, preserves 
fibrations, trivial fibrations, and
weak equivalences between fibrant objects. 
\end{Lem}
\begin{proof}
The $S$-sections functor has a left adjoint, obtained by left 
Kan
extension, that preserves cofibrations and weak
equivalences. See \cite[pg. 60]{Jardine} and
\cite[Cor. 3.16]{Mitchell} for the details.
\end{proof}
\par
We use $\mathcal{L}^2$ to denote the sheafification functor for 
presheaves of sets, simplicial presheaves, and presheaves of 
spectra. Let $i \: \mathbf{ShvSpt} \rightarrow
\mathbf{PreSpt}$ be the inclusion
functor, which is right adjoint to $\mathcal{L}^2$. In \cite[Rk. 
3.11]{GJlocal}, $\mathbf{ShvSpt}$ is given the following model
category structure. A map $h \: \mathcal{F}
\rightarrow \mathcal{G}$ of sheaves of spectra is a weak
equivalence (fibration) 
if and only if $i(f)$ is a weak equivalence (fibration) of 
presheaves.
Also, $h$ is a cofibration
of sheaves of spectra if the following two conditions hold:
\begin{enumerate}
\item 
the map $h^0 \: \mathcal{F}^0 \rightarrow \mathcal{G}^0$ is a 
cofibration of
simplicial presheaves; and
\item
for each $n \geq 0$, the canonical map $\mathcal{L}^2((S^1 \wedge 
\mathcal{G}^n) \cup_{S^1 \wedge \mathcal{F}^n}\mathcal{F}^{n+1})
\rightarrow \mathcal{G}^{n+1}$ is a cofibration of simplicial
presheaves.
\end{enumerate} 
\par
Since $i$ preserves weak equivalences and fibrations between 
sheaves of spectra,
$\mathcal{L}^2$ preserves cofibrations, trivial cofibrations, 
and
weak equivalences between cofibrant objects, and $(\mathcal{L}^2, 
i)$ is a Quillen pair for $(\mathbf{PreSpt},\mathbf{ShvSpt})$. 
By \cite[Cor. 6.22]{Jardine}, if $\mathcal{F}$ is a presheaf of 
sets, a simplicial presheaf, or a presheaf of spectra, then 
\begin{equation}\label{formula}\zig \mathcal{L}^2\mathcal{F} 
\cong \mathrm{Hom}_G(-, \mathrm{colim} \, _N \, \mathcal{F}(G/N)).
\end{equation} This implies that, for any presheaf of spectra 
$\mathcal{F}$, $\mathcal{F} \rightarrow
\mathcal{L}^2 \mathcal{F}$ is a weak equivalence, and thus, 
$\mathrm{Ho}(\mathbf{PreSpt}) \cong \mathrm{Ho}(\mathbf{ShvSpt})$ 
is a Quillen equivalence.
\par
There is an equivalence of categories $\mathbf{ShvSpt} \cong 
Sp_G$, via the functors \[L \: \mathbf{ShvSpt} \rightarrow
Sp_G, \ \ \ L(\mathcal{F}) = \mathrm{colim} \, _N \,
\mathcal{F}(G/N), \ \ \ \mathrm{and}\] \[R \: Sp_G \rightarrow 
\mathbf{ShvSpt}, \ \ \ R(X)= \mathrm{Hom}_G(-, X).\] It is not 
hard to verify this equivalence, since it is an extension to 
spectra of the fact that $\mathbf{Shv}$ and $T_G$ are equivalent 
categories (see \cite[Prop. 6.20]{Jardine}, 
\cite[III-9, Thm. 1]{Maclane}). 
\par
Exploiting this equivalence, we make $Sp_G$ a model category 
in the following
way. Define a map
$f$ of discrete $G$-spectra to be a weak equivalence (fibration) 
if and only if
$\mathrm{Hom}_G(-,f)$ is a weak equivalence (fibration) of 
sheaves of
spectra. Also, define $f$ to be a cofibration if and only if 
$f$
has the left lifting property with respect to all trivial
fibrations. Thus, $f$ is a cofibration if
and only if $\mathrm{Hom}_G(-,f)$ is a cofibration in
$\mathbf{ShvSpt}$. Using this, it is easy
to show that $Sp_G$ is a model category, and there is a Quillen 
equivalence $\mathrm{Ho}(\mathbf{ShvSpt}) \cong
\mathrm{Ho}(Sp_G).$
\par
Formula (\ref{formula}) allows us to define the model category 
structure of $Sp_G$ without reference to sheaves of spectra, in 
the theorem below, extending the model category structure on the 
category $S_G$, given in \cite[Thm. 1.12]{hGal}, to $Sp_G$.
\begin{Thm}\label{finally}
Let $f \: X \rightarrow Y$ be a map in $Sp_G$. Then $f$ is a weak 
equivalence (cofibration) in $Sp_G$ if and only if $f$ is a weak 
equivalence (cofibration) in $Sp$.
\end{Thm}
\begin{proof}
For weak equivalences, the statement is clearly true. Assume that 
$f$ is a cofibration in $Sp_G$. Since $\mathrm{Hom}_G(-,X_0) 
\rightarrow \mathrm{Hom}_G(-,Y_0)$ is a cofibration of simplicial 
presheaves, evaluation at $G/N$ implies that $X_0^N \rightarrow 
Y_0^N$ is a cofibration of simplicial sets. Thus, $X_0 \cong 
\mathrm{colim} \, _N \, X_0^N \rightarrow \mathrm{colim} \, _N 
\, Y_0^N \cong Y_0$ is a cofibration of simplicial sets.
\par
By (\ref{formula}) and the fact that colimits commute with 
pushouts, \[\mathrm{Hom}_G(-,(S^1 \wedge Y_n) \cup_{S^1 \wedge 
X_n}X_{n+1}) \rightarrow \mathrm{Hom}_G(-,Y_{n+1})\] is a 
cofibration of simplicial presheaves, and hence, the map of 
simplicial sets $\mathrm{colim} \, _N \, ((S^1 \wedge Y_n) 
\cup_{S^1 \wedge X_n}X_{n+1})^N \rightarrow Y_{n+1}$ is a 
cofibration. 
\par
Let $W$ be a simplicial pointed discrete $G$-set. Then $S^1 
\wedge W \cong \mathrm{colim} \, _N \, (S^1 \wedge W^N)$, so 
that $S^1 \wedge W$ is also a simplicial pointed discrete 
$G$-set. Since the forgetful functor $U \: T_G \rightarrow 
\mathbf{Sets}$ is a left adjoint, pushouts in $T_G$ are formed 
in $\mathbf{Sets}$, and thus, there is an isomorphism 
\[\mathrm{colim} \, _N \, ((S^1 \wedge Y_n) \cup_{S^1 
\wedge X_n}X_{n+1})^N \cong (S^1 \wedge Y_n) \cup_{S^1 
\wedge X_n}X_{n+1}\] of simplicial discrete $G$-sets. Hence, 
$(S^1 \wedge Y_n) \cup_{S^1 \wedge X_n}X_{n+1} \rightarrow 
Y_{n+1}$ is a cofibration in $\mathcal{S}$, and $f$ is a 
cofibration in $Sp$. 
\par
The converse follows from the fact that if $W \rightarrow Z$ 
is an injection of simplicial discrete $G$-sets, then 
$\mathrm{Hom}_G(-,W) \rightarrow \mathrm{Hom}_G(-,Z)$ is a 
cofibration of simplicial presheaves.
\end{proof}
\begin{Cor}\label{U}
The functors $(U, \mathrm{Map}_c(G, \negthinspace - 
))$ are a Quillen pair for 
$(Sp_G,Sp)$.
\end{Cor}
\par
Let $\mathrm{t} \: Sp \rightarrow Sp_G$ give a
spectrum trivial $G$-action, so that $\mathrm{t}(X)=X$. The 
right adjoint
of $\mathrm{t}$ is the fixed points functor $(-)^G$. Clearly, 
$\mathrm{t}$ preserves all weak equivalences and cofibrations, 
giving the next result.
\begin{Cor}\label{cof}
The functors $(\mathrm{t},(-)^G)$ are a
Quillen pair for $(Sp, Sp_G)$.
\end{Cor}
\par
We conclude this section with a few more useful facts about 
discrete $G$-spectra.
\begin{Lem}\label{fibrant}
If $f \: X \rightarrow Y$ is a fibration in $Sp_G$,
then it is a fibration in $Sp$. In particular, if $X$ is 
fibrant as a discrete $G$-spectrum, then $X$ is
fibrant as a spectrum. 
\end{Lem}
\begin{proof}
Since $\mathrm{Hom}_G(-,f)$ is a fibration of presheaves of
spectra, $\mathrm{Hom}_G(G/N,f)$ is a fibration of 
spectra for
each open normal subgroup $N$. This implies that $\mathrm{colim} 
\, _N \,
\mathrm{Hom}_G(G/N,f)$ is a fibration of spectra. Then the 
lemma follows from factoring $f$ as $X \cong \mathrm{colim} 
\, _N \, X^N \rightarrow
\mathrm{colim} \, _N \, Y^N \cong Y.$
\end{proof}
\par
The next lemma and its corollary show that the homotopy 
groups of a discrete $G$-spectrum are discrete $G$-modules, 
as expected.

\begin{Lem} 
If $X$ is a pointed Kan complex and a simplicial discrete 
$G$-set, then $\pi_n(X)$ is a discrete $G$-module, for all 
$n \geq 2$.
\end{Lem}
\begin{proof}  
Let $f \: S^n \rightarrow X$ be a pointed map. Since $S^n$ 
has only two non-degenerate simplices, the basepoint $\ast$ 
and the fundamental class $\iota_n \in (S^n)_n$, the map $f$ is 
determined by $f(\iota_n)$. The action of $G$ on $\pi_n(X)$ is 
defined as follows: $g \cdot [f]$ is the homotopy class of the 
pointed map $h \: S^n \rightarrow X$ defined by $h(\iota_n) = 
g\cdot f(\iota_n)$. To show 
that the action is continuous, it suffices to show that, for 
any $[f] \in \pi_n(X)$, the stabilizer $G_{[f]}$ of $[f]$ is an 
open subgroup of $G$, and this is elementary.
\end{proof}
\begin{Cor}\label{pie}
If $X$ is a discrete $G$-spectrum, then $\pi_n(X)$ is a discrete 
$G$-module for any integer $n$.
\end{Cor}
\begin{proof}
Let $X \rightarrow X_f$ be a trivial cofibration with $X_f$ 
fibrant, all in $Sp_G$. Then \[\pi_n(X) \cong \pi_n(X_f) \cong 
\mathrm{colim} \, _{m \, \geq \, \mathrm{max}(2-n,0)} \, 
\pi_{m+n}((X_f)_{m})\] is a discrete $G$-module, since each 
$\pi_{m+n}((X_f)_m)$ is one.
\end{proof}  
The following observation says that
certain elementary constructions with discrete $G$-spectra yield 
discrete
$G$-spectra.
\begin{Lem}
Given a profinite group $G \cong \mathrm{lim} \, _N \, G/N$, let 
$\{X_N\}_N$ be a directed system
  of spectra, such that each $X_N$ is a $G/N$-spectrum and the 
maps
  are $G$-equivariant. Then $\mathrm{colim} \, _N \, X_N$ is a
  discrete $G$-spectrum. If $Y$ is a trivial $G$-spectrum, then 
$(\mathrm{colim} \,  _N \, X_N) \wedge Y
\cong \mathrm{colim} \, _N  \, (X_N \wedge Y)$ is a
  $G$-equivariant isomorphism of discrete $G$-spectra. Thus, if 
$X$ is
  a discrete $G$-spectrum, then $X \wedge Y$ is a discrete 
$G$-spectrum.
\end{Lem}
\par 
The corollary below is very useful later.
\begin{Cor}
The spectra $F_n$, $F_n \wedge M_I$, and $F_n \wedge M_I \wedge X$, 
for any spectrum $X$, are discrete $G_n$-spectra.
\end{Cor}
\section{Towers of discrete $G$-spectra and continuous $G$-spectra}
Let $\mathbf{tow}(Sp_G)$ be the category where a
typical object
$\{X_i\}$ is a tower \[\cdots \rightarrow X_i \rightarrow X_{i-1} 
\rightarrow \cdots \rightarrow X_1 \rightarrow X_0\] in $Sp_G$. 
The morphisms are
natural transformations $\{X_i\} \rightarrow \{Y_i\}$, such that
each $X_i \rightarrow Y_i$ is $G$-equivariant. Since $Sp_G$ is a 
simplicial
model category, \cite[VI, Prop. 1.3]{GJ} shows that
$\mathbf{tow}(Sp_G)$ is also a simplicial model category, where
$\{f_i\}$ is a weak equivalence (cofibration) if and only if each
$f_i$ is a weak equivalence (cofibration) in $Sp_G$. By
\cite[VI, Rk. 1.5]{GJ}, if $\{X_i\}$ is fibrant in
$\mathbf{tow}(Sp_G)$, then each map $X_i \rightarrow
X_{i-1}$ in the tower is a fibration and each $X_i$ is fibrant, 
all in
$Sp_G$.       
\begin{Def}\label{fibrep}
Let $\{X_i\}$ be in $\mathbf{tow}(Sp_G)$. Then $\{X'_i\}$ denotes 
the target of a trivial cofibration $\{X_i\} \rightarrow
\{X'_i\},$ with $\{{X'_i}\}$ fibrant, in $\mathbf{tow}(Sp_G).$
\end{Def}
\begin{Rk}
Let $\{X_\alpha\}$ be a diagram in $Sp_G$. Since there is an 
isomorphism $\mathrm{lim}\, ^{Sp_G} _\alpha \, X_\alpha \cong 
\mathrm{colim} \, _N
\, (\mathrm{lim} \, ^{Sp}  _\alpha \, X_\alpha)^N,$ limits in 
$Sp_G$ are not formed in $Sp$. In everything that follows, 
$\mathrm{lim}$ and
$\mathrm{holim}$ are always in $Sp$. 
\end{Rk}
\par
The functor $\mathrm{lim} \, _i \, (-)^G \:
\mathbf{tow}(Sp_G) \rightarrow Sp,$ given by $\{X_i\} \mapsto 
\mathrm{lim} \, _i \, X_i^G,$ is right
adjoint to the functor $\underline{\mathrm{t}} \: Sp \rightarrow
\mathbf{tow}(Sp_G)$ that sends a spectrum $X$ to
the constant diagram $\{X\},$ where $X$ has trivial
$G$-action. Since $\underline{\mathrm{t}}$ preserves all weak
equivalences and cofibrations, we have the following fact.
\begin{Lem}\label{helpful}
The functors $(\underline{\mathrm{t}}, \mathrm{lim} \, _i \, 
(-)^G)$
are a Quillen pair for the categories $(Sp,\mathbf{tow}(Sp_G)).$
\end{Lem} 
\par
This Lemma implies the existence of the total right derived
functor \[\mathbf{R}(\mathrm{lim} \, _i \, (-)^G) \:
\mathrm{Ho}(\mathbf{tow}(Sp_G)) \rightarrow
\mathrm{Ho}(Sp), \ \ \ \{X_i\} \mapsto \mathrm{lim} \,
_i \, ({X'_i})^G.\]
\begin{Lem}\label{ctsspt}
If $\{X_i\}$ in $\mathbf{tow}(Sp_G)$ is a tower of fibrant 
spectra, then there are weak equivalences $\mathrm{holim} \,
  _i \, X_i \overset{p}{\longrightarrow} \mathrm{holim} \,
  _i \, X'_i \overset{q}{\longleftarrow} \mathrm{lim} \, _i \, 
X'_i.$
\end{Lem} 
\begin{proof}
Since each $X_i \rightarrow X'_i$ is a weak equivalence between 
fibrant spectra, $p$ is a weak equivalence. Since $X'_i
\rightarrow X'_{i-1}$
is a fibration in $Sp$, for $i \geq 1$, $q$ is a
weak equivalence.
\end{proof}
\par
We use this lemma to define the continuous $G$-spectra that we 
will study. 
\begin{Def}\label{ctsdef}
If $\{X_i\}$ is in $\mathbf{tow}(Sp_G)$, then the
inverse limit $\mathrm{lim} \, _i \, X_i$ is a {\em continuous 
$G$-spectrum}. Also, if $\{X_i\} \in \mathbf{tow}(Sp_G)$ is a 
tower of fibrant spectra, we call $\mathrm{holim} \, _i \, X_i$ a 
{\em continuous $G$-spectrum}, due to the zigzag of Lemma 
\ref{ctsspt} that relates $\mathrm{holim} \, _i \, X_i$ to the 
continuous $G$-spectrum $\mathrm{lim} \, _i \, X'_i.$ 
Notice that if $X \in Sp_G$, then, using the constant tower on $X$, 
$X \cong \mathrm{lim} \, _i \, X$ is a continuous $G$-spectrum.
\end{Def}
\begin{Rk}\label{loose}
Sometimes we will use the term
``continuous $G$-spectrum'' more loosely. Let $X$ be a continuous
$G$-spectrum, as in Definition \ref{ctsdef}. If $Y$ is a 
$G$-spectrum 
that is isomorphic to $X$, in
the stable category, with compatible $G$-actions, then we 
call $Y$ a continuous $G$-spectrum.
\end{Rk} 
\par
We make a few comments about Definition \ref{ctsdef}. Though the 
definition is not as general as it could be, it is sufficient for 
our applications. The inverse limit is central to the definition 
since the inverse limit of a tower of discrete $G$-sets is a 
topological $G$-space. 
\par
Given any tower $X_\ast = \{X_i\}$ in $Sp_G$, \[\mathrm{holim} \, 
_i \, X_i = \mathrm{Tot}(\textstyle{\prod}^\ast X_\ast) \cong 
\mathrm{lim} \, _n \, \mathrm{T}(n),\] where $\mathrm{T}(n) = 
\mathrm{Tot}_n(\prod^\ast X_\ast).$ (See \S \ref{disdef} for the 
definition of $\prod^\ast X_\ast$, and \cite{Bousfield/Kan} is a 
reference for any undefined notation in this paragraph.) Then it 
is natural to ask if $T(n)$ is in $Sp_G$, so that $\mathrm{holim} 
\, _i \, X_i$ is canonically a continuous $G$-spectrum. For this 
to be true, it must be that, for any $m \geq 0$, the simplicial 
set $\mathrm{T}(n)_m = \mathrm{Map}_{c\mathcal{S}}(\mathrm{sk}_n 
\Delta[-], \prod^\ast (X_\ast)_m)$ belongs to $\mathcal{S}_G$. 
If, for all $s \geq 0$, $\prod^s (X_\ast)_m \in \mathcal{S}_G$, 
then $\mathrm{T}(n)_m \in \mathcal{S}_G$, by \cite[pg. 212]{hGal}. 
However, the infinite product $\prod^s (X_\ast)_m$ need not be in 
$\mathcal{S}_G$, and thus, in general, $\mathrm{T}(n) \notin 
Sp_G$. Therefore, $\mathrm{holim} \, _i \, X_i$ is not always 
identifiable with a continuous $G$-spectrum, in the above way.
\section{Homotopy fixed points for discrete $G$-spectra}
\label{disdef}
In this section, we define the homotopy fixed point spectrum 
for
$X \in Sp_G$. We begin by recalling the homotopy spectral 
sequence, since we use it often.
\par
If $J$ is a small category, $P \: J
\rightarrow Sp$ a diagram of fibrant spectra, and $Z$ any 
spectrum, then there is a conditionally convergent spectral
sequence \begin{equation}\label{intobkss}\zig
E_2^{s,t} = \mathrm{lim}^s \,
_J \, [Z,P]_{t} \Rightarrow [Z, \mathrm{holim} \, _J \,
P]_{t-s},\end{equation} where $\mathrm{lim}^s$ is the $s$th right 
derived functor of $\mathrm{lim} \,
_J  \: \mathbf{Ab}^J \rightarrow \mathbf{Ab}$ \cite[Prop. 5.13, 
Lem.
5.31]{Thomason}.  
\par
Associated to $P$ is
the cosimplicial spectrum $\prod^\ast P$, with 
$\textstyle{\prod}^n P =
\textstyle{\prod} \, _{(B\Delta)_n} \, P({j_n}),$ where the
$n$-simplices of the nerve $B\Delta$ consist of all
strings $[j_0] \rightarrow \cdots \rightarrow [j_n]$ of $n$
morphisms in $\Delta$. For any $k \geq 0$,
$(\textstyle{\prod}^\ast P)_k$ is a fibrant cosimplicial
pointed simplicial set, and $\mathrm{holim} \, _J \,
P = \mathrm{Tot}(\textstyle{\prod}^\ast P).$ 
\par
If $P$ is a cosimplicial diagram of fibrant spectra (a
{\em cosimplicial fibrant spectrum}), then $E_2^{s,t} = \pi^s 
[Z,P]_t,$ the $s$th cohomotopy group of the cosimplicial abelian 
group $[Z,P]_t$.
\begin{Def}\label{profinite}
Given a discrete $G$-spectrum $X$, $X^{hG}$ denotes the {\em 
homotopy fixed
point spectrum\/} of $X$ with respect to the continuous action of
$G$. We define $X^{hG} = (X_{f,G})^G,$ where $X \rightarrow 
X_{f,G}$ is a
trivial cofibration and $X_{f,G}$ is fibrant, in $Sp_G$. We 
write $X_f$ instead of $X_{f,G}$ when doing so causes no 
confusion.
\end{Def}
\par
Note that $X^{hG} = \mathrm{Hom}_G(\ast,X_{f,G})$, the global 
sections of the presheaf of spectra $\mathrm{Hom}_G(-,X_{f,G}).$ 
This definition has been developed in other categories: see 
\cite{hGal} for simplicial discrete
$G$-sets, \cite{Jardinejpaa} for simplicial presheaves, and 
\cite[Ch. 6]{Jardine} for presheaves of spectra.
\par
As expected, the definition of homotopy fixed points for a
profinite group generalizes the definition for a finite
group. Let $G$ be a finite group, $X$ a
$G$-spectrum, and let $X \rightarrow X_\mathtt{f}$ be a weak 
equivalence that
is $G$-equivariant, with $X_\mathtt{f}$ a fibrant spectrum. 
Then, since $G$ is finite, the homotopy fixed point spectrum 
$X^{h'G} = \mathrm{Map}_G(EG_+,X_\mathtt{f})$ can be defined to 
be $\mathrm{holim} \, _G \, X_\mathtt{f}$ (as in the Introduction). 
Note that there is a descent spectral sequence \[E_2^{s,t} = 
\textstyle{\mathrm{lim}}^s \, _G \, \pi_t(X) \cong
H^s(G;\pi_t(X)) \Rightarrow \pi_{t-s}(X^{h'G}).\] Since $G$ is
profinite, $X$ is a discrete $G$-spectrum, and $X_{f,G}$ is a
fibrant spectrum, so that $X^{h'G} = \mathrm{holim} \, _G \,
X_f.$ Then, by \cite[Prop. 6.39]{Jardine}, the canonical map 
$X^{hG} =
\mathrm{lim} \, _G \, X_f \rightarrow \mathrm{holim} \, _G \, X_f
= X^{h'G}$ is a weak equivalence, as desired. 
\par
We point out several properties of homotopy fixed points that 
follow from Corollary \ref{cof}.
\begin{Lem}\label{tright}
The homotopy fixed points functor $(-)^{hG} \: \mathrm{Ho}(Sp_G)
\rightarrow \mathrm{Ho}(Sp)$ is the total right derived functor of
the fixed points functor $(-)^G \: Sp_G \rightarrow
Sp.$ In particular, if $X \rightarrow Y$ is a weak equivalence 
of discrete
$G$-spectra, then $X^{hG} \rightarrow
Y^{hG}$ is a weak equivalence.
\end{Lem}
\par
Given two globally fibrant models $\mathbf{G}\mathcal{F}$ and
$\mathbf{G'}\mathcal{F}$ for a presheaf of spectra
$\mathcal{F}$, there need not be a weak equivalence
$\mathbf{G}\mathcal{F} \rightarrow \mathbf{G'}\mathcal{F},$ 
since not every presheaf of spectra is cofibrant. Also,
a cofibration of sheaves of spectra can fail to be a cofibration 
of
presheaves. Given these limitations, there is the following
relationship between homotopy fixed points and the global
sections of globally fibrant models, due to the left lifting 
property of trivial cofibrations with respect to fibrations. 
\begin{Lem}
Let $X \in Sp_G$. Then $X^{hG} \overset{\simeq}{\longleftarrow} 
\mathcal{F}(\ast)
\overset{\simeq}{\longrightarrow}
\mathbf{G}\mathrm{Hom}_G(-,X)(\ast),$ where $\mathrm{Hom}_G(-,X) 
\rightarrow \mathcal{F}$ is
a trivial cofibration of presheaves, with $\mathcal{F}$ globally 
fibrant. Thus, $X^{hG}$ and
$\mathbf{G}\mathrm{Hom}_G(-,X)(\ast)$ have the same stable 
homotopy
type.
\end{Lem}
\section{$E_n$ is a continuous $G_n$-spectrum}\label{7}
We show that $E_n$ is a continuous $G_n$-spectrum by successively
eliminating simpler ways of constructing a continuous
action, and by applying the theory of the previous section.
\par
First of all, since the profinite ring $\pi_0(E_n)$ is not a 
discrete $G_n$-module, Corollary \ref{pie} implies the following 
observation.
\begin{Lem}
$E_n$ is not a discrete $G_n$-spectrum.
\end{Lem}
However, note that, for $k \in \Zset$,
$\pi_{2k}(E_n \wedge M_I)$ is a finite discrete $G_n$-module so 
that the
action factors through a finite quotient
$G_n/{U_I}$, where $U_I$ is some
open normal subgroup (see \cite[Lem. 1.1.16]{Ribes}). Thus, 
$\pi_{2k}(E_n \wedge M_I)$ is a $G_n/{U_I}$-module, and one is 
led
to ask if $E_n \wedge M_I$ is a $G_n/{U_I}$-spectrum. If
so, then $E_n \wedge M_I$ is a discrete $G_n$-spectrum, and $E_n$
is easily seen to be a continuous $G_n$-spectrum. 
\par
However, $E_n \wedge M_I$ is not a $G_n/U$-spectrum for all open 
normal
subgroups $U$ of $G_n$, as the following observation shows. As 
far as
the author knows, this observation is due to Mike Hopkins; the
author learned the proof from Hal Sadofsky.
\begin{Lem}\label{Sadofsky} There is no open normal subgroup 
$U$ of $\GG$ such that
  the $\GG$-action on $E_n \wedge M_I$ factors through $\GG/U$.
\end{Lem}
\begin{proof}
Suppose the $\GG$-spectrum $E_n \wedge M_I$ is a
$\GG/U$-spectrum. Then the $\GG$-action on the
middle factor of $E_n \wedge E_n \wedge M_I$ factors through
$\GG/U$, so that $\pi_\ast(E_n \wedge E_n \wedge M_I)$ is a 
$\GG/U$-module. Note that $\pi_\ast(E_n \wedge E_n \wedge
M_I) \cong \mathrm{Map}^\ell_c(\GG,E_{n \ast}/I)$.
\par
Since the $\GG$-module $\mathrm{Map}^\ell_c(\GG,E_{n \ast}/I)$ 
is
a $\GG/U$-module, there is an isomorphism of sets
$\mathrm{Map}^\ell_c(\GG,E_{n,0}/I) =
\mathrm{Map}^\ell_c(\GG,E_{n,0}/I)^U \cong \mathrm{Map}_c(G_n/U,
E_{n,0}/I).$ But the first set is infinite and the last is
finite, a contradiction. 
\end{proof}
\par
Since $\pi_\ast(E_n \wedge
M_I)$ is a discrete $\GG$-module, one can still hope for a 
spectrum
$E_n/I \simeq E_n \wedge M_I$, such that $E_n/I$ is a discrete 
$\GG$-spectrum. 
\par
To produce $E_n/I$, we make the following observation. 
By \cite[Remark 6.26]{Jardine}, since $U_i$ is an open 
normal subgroup of $G_n$, the presheaf 
$\mathrm{Hom}_{U_i}(-,(E_n/I)_{f,G_n})$ is fibrant in 
the model category of presheaves of spectra on the site 
$U_i \negthinspace - \negthinspace \mathbf{Sets}_{df}$. 
Thus, for each $i$, the map $E_n/I \rightarrow (E_n/I)_{f,G_n}$ 
is a trivial cofibration, with fibrant target, all in $Sp_{U_i}$, 
so that $((E_n/I)_{f,G_n})^{U_i} = (E_n/I)^{hU_i}.$
\par
Combining this observation with the idea,
discussed in \S 1, that $E_n \wedge M_I$ has homotopy fixed 
point spectra $(E_n \wedge M_I)^{hU_i} \simeq E_n^{\mathtt{h}U_i} 
\wedge M_I \simeq (E_n/I)^{hU_i},$ we have:
\begin{align*}E_n/I & \simeq (E_n/I)_{f,G_n} \cong \mathrm{colim} 
\, _i \, ((E_n/I)_{f,G_n})^{U_i} = \mathrm{colim} \, _i \,
(E_n/I)^{hU_i} \\ & \simeq \mathrm{colim} \, _i \,
(E_n^{\mathtt{h}U_i} \wedge M_I) \cong F_n \wedge M_I.\end{align*} 
This argument suggests that $E_n \wedge M_I$ has the homotopy 
type of the discrete $G_n$-spectrum $F_n \wedge M_I$. To show 
that this is indeed the case, we consider the spectrum $F_n$ 
in more detail. The key result is the following theorem, 
due to Devinatz and Hopkins.
\par 
\begin{Thm}[{\cite{DH}}]\label{central} 
There is a weak equivalence $E_n \simeq \Lhat (F_n).$
\end{Thm}
\begin{proof}
By \cite[Thm. 3]{DH}, $E_n^{h'\{e\}}
\simeq E_n^{\mathtt{h}\{e\}}.$ (We remark that this 
weak equivalence is far from obvious.) By 
\cite[Definition 1.5]{DH}, $E_n^{\mathtt{h}\{e\}} = \Lhat 
(\mathrm{hocolim} \,  _i \, E_n^{\mathtt{h}U_i}),$ 
where the homotopy colimit is in the category $E_\infty$. 
Then, by \cite[Remark 1.6, Lemma 6.2]{DH},  
$\mathrm{hocolim} \,  _i \, E_n^{\mathtt{h}U_i} 
\simeq \mathrm{colim} \, _i \, E_n^{\mathtt{h}U_i},$ 
where the colimit is in $\mathcal{M}_S.$ Thus, as 
spectra in $Sp$, $E_n^{\mathtt{h}\{e\}} \simeq \Lhat 
(F_n)$, so that $E_n \simeq E_n^{h'\{e\}} \simeq \Lhat (F_n).$
\end{proof}
\begin{Cor}\label{key}
In the stable category, there are isomorphisms \[E_n
\cong \mathrm{holim} \, _I \, (F_n \wedge M_I)_\mathtt{f} 
\cong
\mathrm{holim} \, _I \, (E_n \wedge M_I)_\mathtt{f}.\]
\end{Cor}
\par
The following result, a proof of which was first shown to me by 
Charles Rezk, shows that $E_n \wedge M_I \simeq F_n
\wedge M_I$. This weak equivalence and $\mathrm{vcd}(G_n) < 
\infty$ 
are the main facts that make it possible to construct the homotopy 
fixed point spectra
of $E_n$.
\begin{Cor}\label{Rezk} 
If $Y$ is a finite spectrum
  of type $n$, then the $G_n$-equivariant map $F_n \wedge Y
  \rightarrow E_n \wedge Y$ is a weak equivalence. In 
particular,
  $E_n \wedge M_I \simeq F_n \wedge M_I$.
\end{Cor}
\begin{proof} We have $E_n \wedge Y \simeq \Lhat  (F_n)
  \wedge Y \simeq L_n(F_n) \wedge Y \simeq F_n \wedge Y.$
\end{proof}
\par
Now we show that $E_n$ is a continuous $G_n$-spectrum.
\begin{Thm}\label{wonder}
There is an isomorphism $E_n \cong \mathrm{holim} \, _I \, (F_n
\wedge M_I)_{f,G_n}$. Thus, $E_n$ is a continuous $G_n$-spectrum.
\end{Thm}
\begin{proof}
By Corollary \ref{key}, $E_n \cong 
\mathrm{holim} \, _I \, (F_n \wedge M_I)_\mathtt{f}$. By 
functorial fibrant replacement, the 
map of towers $\{F_n \wedge M_I\} \rightarrow \{(F_n 
\wedge M_I)_{f,G_n}\}$ induces a map of towers $\{(F_n \wedge 
M_I)_\mathtt{f}\} \rightarrow 
\{((F_n \wedge M_I)_{f,G_n})_\mathtt{f}\}$ and, hence, 
weak equivalences 
\[\mathrm{holim} \, _I \, (F_n
\wedge M_I)_{f,G_n} \rightarrow \mathrm{holim} \, _I \, 
((F_n
\wedge M_I)_{f,G_n})_\mathtt{f} \leftarrow \mathrm{holim} 
\, _I \, (F_n
\wedge M_I)_\mathtt{f}.\] Thus, $\mathrm{holim} \, _I \, 
(F_n \wedge M_I)_{f,G_n}$ is isomorphic to $\mathrm{holim} 
\, _I \, (F_n
\wedge M_I)_\mathtt{f}$ and $E_n$. Since 
$\{(F_n \wedge M_I)_{f,G_n}\}$ is a tower of fibrant 
spectra, $\mathrm{holim} \, _I \, (F_n
\wedge M_I)_{f,G_n}$ is a continuous $G_n$-spectrum. 
Then, by Remark \ref{loose}, $E_n$ is a continuous 
$G_n$-spectrum.
\end{proof}
We conclude this section with some observations about $F_n$.
\begin{Lem}
The map $\eta \: F_n \rightarrow E_n$ is not a weak equivalence 
and $F_n$ is not $K(n)$-local.
\end{Lem}
\begin{proof}
If $\eta$ is a weak equivalence, then $\pi_0(\eta)$ is a 
$\GG$-equivariant
isomorphism from a discrete $G_n$-module (with all orbits 
finite) to a non-finite profinite $G_n$-module, which
is impossible. If $F_n$ is $K(n)$-local, then $F_n \simeq \Lhat  
(F_n) \simeq E_n$, and $\eta$ is a weak equivalence, a 
contradiction.
\end{proof}
\begin{Lem}
The maps $\Lhat  (F_n
\wedge F_n) \rightarrow  \Lhat  (E_n
\wedge F_n) \rightarrow \Lhat  (E_n
\wedge E_n)$ are weak equivalences.
\end{Lem}
\begin{proof}
Since $F_n \wedge M_I \simeq E_n \wedge M_I$, $F_n \wedge
F_n \wedge M_I \simeq E_n \wedge E_n \wedge M_I.$
Since $F_n \wedge F_n$, $E_n \wedge F_n$ and $E_n \wedge E_n$ 
are
$E(n)$-local, the result follows from Theorem \ref{veryneat}.
\end{proof}
\par 
Since $E_n^{\mathtt{h}U_i}$ is not an $E_n$-module (but $E_n$ 
is an $E_n^{\mathtt{h}U_i}$-module), note that $F_n$ and 
$\mathrm{holim} \, _I \, (F_n \wedge M_I)_{f,G_n}$ are not 
$E_n$-modules. However, any fully satisfactory model of $E_n$ 
as a continuous $\GG$-spectrum would
also be a twisted $E_n$-module $\GG$-spectrum, since $E_n$ 
itself is such a spectrum. (``Twisted'' means that the module 
structure map is
$\GG$-equivariant, where $E_n \wedge E_n$ has the diagonal 
action.)
\section{Homotopy fixed points when $\mathrm{vcd(G)} \, 
\mathrm{<} \, \infty$}\label{hfps}
In this section, $G$ always has finite virtual cohomological 
dimension. Thus, there exists a uniform bound $m$, such that 
for all $K <_c G$, $\mathrm{vcd}(K) \leq m.$  For $X \in Sp_G$, 
we use this fact to give a model for $X^{hG}$ that eases the 
construction of its descent spectral sequence. Also, this fact 
yields a second model for $X^{hK}$ that is functorial in $K$. 
\begin{Def} Consider the functor 
\[\Gamma_G = \mathrm{Map}_c(G,-) \circ U \: Sp_G \rightarrow
Sp_G, \ \ \ X \mapsto \Gamma_G(X)=
\mathrm{Map}_c(G,X),\] where $\Gamma_G(X)$ has the $G$-action
defined in \S\ref{spt}. We write $\Gamma$ instead of $\Gamma_G$,
when $G$ is understood from context. There is a
$G$-equivariant monomorphism $i \: X \rightarrow \Gamma
X$ defined, on the level of sets, by $i(x)(g)= g \cdot x.$ 
As in \cite[8.6.2]{Weibel}, since $U$ and $\mathrm{Map}_c(G,-)$ 
are adjoints, $\Gamma$ forms a triple and there is a 
cosimplicial discrete $G$-spectrum
$\Gamma^\bullet X,$ with $(\Gamma^\bullet X)^k \cong 
\mathrm{Map}_c(G^{k+1},X)$.
\end{Def}
\par
We recall the construction of Thomason's hypercohomology
spectrum for the topos $\mathbf{Shv}$ 
(see \cite[1.31-1.33]{Thomason} and 
\cite[\S 1.3, \S 3.2]{Mitchell} for more details). 
Consider the functor 
\[T=u_\ast u^\ast \: \mathbf{ShvSpt}
\rightarrow \mathbf{ShvSpt}, \ \ \ \negthinspace 
\mathcal{F} \mapsto 
\mathrm{Hom}_G(-, \mathrm{Map}_c(G,
\mathrm{colim} \, _N \, \mathcal{F}(G/N))),\] obtained by 
composing the adjoints in the point of the topos. Then, for 
$X \in Sp_G$, $T (\mathrm{Hom}_G(-,X)) \cong
\mathrm{Hom}_G(-,\mathrm{Map}_c(G, X)).$ By iterating this 
isomorphism, the cosimplicial sheaf of spectra $T^\bullet 
\mathrm{Hom}_G(-,X)$ gives rise to the cosimplicial sheaf 
$\mathrm{Hom}_G(-, \Gamma^\bullet X)$. 
\begin{Def}
Given $X \in Sp_G$, the 
{\em presheaf of hypercohomology spectra of
$G \negthinspace - \negthinspace \mathbf{Sets}_{df}$ with 
coefficients in $X$\/} is the presheaf of spectra 
\[\mathbb{H}^\bullet \negthinspace (-;X) = \mathrm{holim} \, 
_\Delta \, \mathrm{Hom}_G(-, \Gamma^\bullet X) \: 
(G\negthinspace - \negthinspace 
\mathbf{Sets}_{df})^{\mathrm{op}}
\rightarrow Sp,\] 
and $\mathbb{H}^\bullet \negthinspace (S;X) = \mathrm{holim} \, 
_\Delta \, \mathrm{Hom}_G(S, \Gamma^\bullet X)$ is the 
hypercohomology spectrum of $S$ with
coefficients in $X$. 
\end{Def}
\par
The map $X_f \rightarrow \Gamma^\bullet
X_f$, induced by $i$, out of the constant cosimplicial diagram, 
and $\mathrm{Hom}_G(-, X_f) \rightarrow \mathrm{lim} \, _\Delta 
\, \mathrm{Hom}_G(-, X_f) \rightarrow
 \mathrm{holim} \, _\Delta \, \mathrm{Hom}_G(-, X_f)$ induce a 
canonical map $\mathrm{Hom}_G(-, X_f) \rightarrow 
\mathbb{H}^\bullet \negthinspace (-,X_f).$
\par
Now we show that $\mathbb{H}^\bullet \negthinspace(\ast,X_f)$ 
is a model for $X^{hG}$. Below, a {\em cosimplicial globally 
fibrant presheaf\/} is a cosimplicial presheaf of spectra that is 
globally fibrant at each level.
\begin{Lem}\label{globally}
If $\mathcal{F}^\bullet$ is a cosimplicial globally fibrant 
presheaf, then $\mathrm{holim} \, _\Delta \, \mathcal{F}^\bullet$ 
is a globally fibrant presheaf. 
\end{Lem}
\begin{proof}
By \cite[Rk. 2.35]{Jardine}, this is equivalent to showing that, 
for each $n \geq 0$, (a) $\mathrm{holim} \, _\Delta \, 
(\mathcal{F}^\bullet)^n$ is a globally fibrant simplicial 
presheaf; and (b) the adjoint of the bonding map, the composition 
\[\gamma\: \mathrm{holim} \, _\Delta \, (\mathcal{F}^\bullet)^n 
\rightarrow \Omega(\mathrm{holim} \, _\Delta \, 
(\mathcal{F}^\bullet)^{n+1}) \cong \mathrm{holim} \, _\Delta \, 
\Omega(\mathcal{F}^\bullet)^{n+1}\] is a local weak equivalence 
of simplicial presheaves. Part (a), the difficult part of this 
lemma, is proven in \cite[Prop. 3.3]{Jardinejpaa}.
\par
The map $\gamma$ is a local weak equivalence if the map of stalks 
$\mathrm{colim} \, _N \, \gamma(G/N)$ is a weak equivalence in 
$\mathcal{S}$, which is true if each $\gamma(G/N)$ is a weak 
equivalence. By \cite[pg. 74]{Jardinejpaa}, if $P$ is a globally 
fibrant simplicial presheaf, then $\Omega P$ is too, so that 
$\sigma \: (\mathcal{F}^\bullet)^n  \rightarrow 
\Omega(\mathcal{F}^\bullet)^{n+1}$ is a cosimplicial diagram 
of local weak equivalences between globally fibrant simplicial 
presheaves. Thus, $\sigma(G/N)$ is a cosimplicial diagram of 
weak equivalences between Kan complexes, so that 
$\gamma(G/N)=\mathrm{holim} \, _\Delta \, \sigma(G/N)$ is 
indeed a weak equivalence.
\end{proof}
\par
The following result is not original: it is basically a special 
case of {\cite[Prop. 3.3]{Jardinejpaa}, versions of which appear 
in \cite[\S 5]{hGal}, \cite[Prop. 3.20]{Mitchell}, and 
\cite[Prop. 6.1]{Mitchell2}}. Since the result is central to our 
work, for the benefit of the reader, we give the details of the 
proof. For use now and later, we recall that, for any group $G$ 
and any closed subgroup $K$, $H^s_c(K;\mathrm{Map}_c(G,A)) = 0$, 
when $s>0$ and $A$ is any discrete abelian group 
\cite[Lemma 9.4.5]{Wilson}. 
\begin{Thm}\label{fibrantmodel}
Let $G$ be a profinite group with $\mathrm{vcd}(G) \leq m$,
and let $X$ be a discrete
$G$-spectrum. Then there are weak equivalences
\[\mathrm{Hom}_G(-,X) \overset{\simeq}{\longrightarrow} 
\mathrm{Hom}_G(-,X_f)
\overset{\simeq}{\longrightarrow} \mathrm{holim} \, _\Delta \, 
\mathrm{Hom}_G(-, \Gamma^\bullet X_f),\] and $\mathrm{holim} 
\, _\Delta \, \mathrm{Hom}_G(-, \Gamma^\bullet X_f)$ is a 
globally fibrant model for
$\mathrm{Hom}_G(-,X)$. Thus, evaluation at $\ast \in 
G\negthinspace - \negthinspace \mathbf{Sets}_{df}$ gives a weak 
equivalence $X^{hG} \rightarrow \mathrm{holim}
\, _\Delta \, (\Gamma^\bullet X_f)^G.$
\end{Thm}
\begin{Rk}
The weak equivalence $X \cong \mathrm{colim} \, _N \, X^N 
\rightarrow \mathrm{colim} \, _N \, (X^N)_\mathtt{f},$ whose 
target is a discrete $G$-spectrum that is fibrant in $Sp$, 
induces a weak equivalence $X_{f,G} \rightarrow (\mathrm{colim} 
\, _N \, (X^N)_\mathtt{f})_{f,G}.$ Thus, there are weak 
equivalences 
\[\mathbb{H}^\bullet \negthinspace (\ast,X_f) \rightarrow 
\mathbb{H}^\bullet \negthinspace (\ast, (\mathrm{colim} \, _N 
\, (X^N)_\mathtt{f})_{f,G}) \leftarrow \mathbb{H}^\bullet 
\negthinspace (\ast, \mathrm{colim} \, _N \, (X^N)_\mathtt{f}),\] 
so that $\mathbb{H}^\bullet \negthinspace (\ast, \mathrm{colim} 
\, _N \, (X^N)_\mathtt{f})$ is a model for $X^{hG}$ that does 
not require the model category $Sp_G$ for its construction.
\end{Rk}
\begin{proof}[Proof of Theorem \ref{fibrantmodel}.]
Since $X_f$ is fibrant in $Sp$, $\Gamma X_f$ is fibrant in $Sp_G$ 
by Corollary \ref{U}. By iteration, $\mathrm{Hom}_G(-, 
\Gamma^\bullet X_f)$ is a cosimplicial globally fibrant presheaf, 
so that $\mathrm{holim} \, _\Delta \, \mathrm{Hom}_G(-, 
\Gamma^\bullet X_f)$ is globally fibrant, by Lemma \ref{globally}. 
It
only remains to show that $\lambda \: \mathrm{Hom}_G(-,X) 
\rightarrow \mathrm{holim} \, _\Delta \, \mathrm{Hom}_G(-, 
\Gamma^\bullet X_f)$ is a weak
equivalence. 
\par
By hypothesis, $G$ contains an open subgroup $H$ with 
$\mathrm{cd}(H) \leq m$. Then by \cite[Lem. 0.3.2]{Wilson}, 
$H$ contains a subgroup
$K$ that is an open normal subgroup of $G$. Let $\{N\}$ be the
collection of open normal subgroups of $G$. Let $N' = N \cap K$. 
Observe that $\{N'\}$ is a cofinal subcollection of open normal 
subgroups of
$G$ so that $G \cong \mathrm{lim} \, _{N'} \, G/{N'}.$ Since 
$N'<_c H$, $\mathrm{cd}(N')
\leq \mathrm{cd}(H).$ Thus, $H^s_c(N'; M) = 0,$ for all
$s>m+1$, whenever $M$ is a discrete $N'$-module. Henceforth, 
we drop the $'$ from $N'$ to ease the notation: $N$
is really $N \cap K$.
\par
Any presheaf of sets $\mathcal{F}$ has stalk $\mathrm{colim}
\, _{N} \, \mathcal{F}(G/N)$, so that $\lambda$ is a weak
equivalence if $\lambda_u \: X \cong
\mathrm{colim}
\, _{N} \, X^{N} \rightarrow
\mathrm{colim} \, _{N} \, \mathrm{holim} \, _\Delta \,
(\Gamma^\bullet X_f)^{N}$ is a weak equivalence. Since
$\mathrm{Hom}_G(-,
\Gamma^\bullet X_f)$ is a cosimplicial globally fibrant 
spectrum, 
the diagram $(\Gamma^\bullet X_f )^{N}$ is a cosimplicial 
fibrant
spectrum. Then, for each $N$, there is a conditionally
convergent spectral
sequence \begin{equation}\label{specseq}\zig
E_2^{s,t}(N) = \pi^s \pi_t ((\Gamma^\bullet X_f )^{N})
\Rightarrow \pi_{t-s}(\mathrm{holim} \, _\Delta \, 
(\Gamma^\bullet X_f)^{N}).\end{equation} 
Because $\pi_t(X)$ is a discrete $G$-module, we have
\[\pi_t(\mathrm{Map}_c(G,
X_f)^{N}) \cong \pi_t(\mathrm{Map}_c(G/{N}, X_f))
\cong \textstyle{\prod} \, _{G/{N}} \, \pi_t(X) \cong 
\mathrm{Map}_c(G,
\pi_t(X))^{N}\] and \[\pi_t(\mathrm{Map}_c(G, X_f)) \cong 
\pi_t(\mathrm{colim} \, _N \, \textstyle{\prod} \, _{G/N} \, 
X_f) \cong \mathrm{Map}_c(G,\pi_t(X)).\] Iterating such 
manipulations, we obtain $\pi^s \pi_t
((\Gamma^\bullet X_f )^{N}) \cong
H^s((\Gamma^\ast \pi_t(X))^{N}).$ The
cochain complex $0 \rightarrow \pi_t(X) \rightarrow 
\Gamma^\ast \pi_t(X)$ of discrete $N$-modules is exact 
(see e.g. 
\cite[pp. 210-211]{hGal}), and, for $k \geq 1$ and $s >0,$ 
\[H^s_c(N; \Gamma^k \pi_t(X)) \cong H^s_c(N;
\mathrm{Map}_c(G, \Gamma^{k-1} \pi_t(X))) 
= 0.\] Thus, the above cochain complex is a resolution
of $\pi_t(X)$ by $(-)^{N}$-acyclic modules, so that 
$E_2^{s,t}(N)
\cong H^s_c(N;
\pi_t(X)).$
Taking a colimit over $\{N \}$ of
(\ref{specseq}) gives the spectral sequence
\begin{equation}\label{specseq2}\zig
E_2^{s,t} = \mathrm{colim} \, _{N} \, H^s_c(N; \pi_t(X))
\Rightarrow \pi_{t-s}(\mathrm{colim} \, _{N} \, \mathrm{holim} 
\, _\Delta \, (\Gamma^\bullet X_f)^{N}).\end{equation} Since
$E_2^{s,\ast}(N) = 0$ whenever $s>m+1$, the $E_2$-terms
$E_2(N)$ are uniformly
bounded on the right. Therefore, by \cite[Prop. 3.3]{Mitchell},
the colimit of the spectral sequences does converge to the
colimit of the abutments, as asserted in (\ref{specseq2}). 
\par 
Finally, $E_2^{*,t} \cong H^\ast_c(\mathrm{lim} \, _N \, N; 
\pi_t(X)) \negthinspace \cong \negthinspace H^\ast(\{e\};
\pi_t(X)),$ which is isomorphic to $\pi_t(X)$, concentrated 
in
degree zero.
Thus,
(\ref{specseq2}) collapses and for all $t$, 
\[\pi_{t}(\mathrm{colim} \, _{N} \, \mathrm{holim} \, 
_\Delta \,(\Gamma^\bullet X_f)^{N}) \cong
\pi_t(X),\] and hence, $\lambda_u$ is a weak equivalence.
\end{proof}
\begin{Rk}\label{model}
Because of Theorem \ref{fibrantmodel}, if $\mathrm{vcd}(G) < 
\infty$ and $X$ is a discrete $G$-spectrum, we make the 
identification
\[X^{hG} = \mathrm{holim} \, _\Delta \, (\Gamma^\bullet_G 
X_{f,G})^G = \mathbb{H}^\bullet \negthinspace (\ast,X_{f,G}).\]
\end{Rk}
\par
In \cite[\S 2.14]{Stavros}, an expression that is basically 
equivalent to $\mathrm{holim} \, _\Delta \, 
(\Gamma^\bullet_G X_{f,G})^G$ is defined to be the homotopy 
fixed point spectrum $X^{hG}$, even if $\mathrm{vcd}(G) = 
\infty.$ This approach has the disadvantage that $(-)^{hG}$ 
need not always be the total right derived functor of 
$(-)^G$. Thus, we only make the identification of 
Remark \ref{model} when $\mathrm{vcd}(G) < \infty.$
\par 
Now it is easy to construct the descent spectral sequence. 
Note 
that if $X$ is a discrete $G$-spectrum, the proof of Theorem 
\ref{fibrantmodel} 
shows that \[\pi^s \pi_t ((\Gamma^\bullet_G X)^G) \cong 
\pi^s((\Gamma^\bullet_G \pi_t(X))^G) \cong H^s_c(G;\pi_t(X)).\]
\begin{Thm}\label{vcddss}
If $\, \mathrm{vcd}(G) < \infty$ and $X$ is a discrete 
$G$-spectrum, then there is a
conditionally convergent descent spectral sequence
\[E_2^{s,t} = H^s_c(G; \pi_t(X)) \Rightarrow 
\pi_{t-s}(X^{hG}).\]
\end{Thm}
\begin{proof} 
As in Theorem \ref{fibrantmodel}, $(\Gamma^\bullet X_f)^G$ is a 
cosimplicial fibrant spectrum. Thus, we can form the homotopy 
spectral sequence for $\pi_\ast(\mathrm{holim} \, _\Delta \,  
(\Gamma^\bullet X_f)^G).$
\end{proof}
\begin{Rk}
The descent spectral sequence of Theorem \ref{vcddss} has 
been constructed in other contexts: for simplicial presheaves, 
presheaves of spectra, and $S_G$, see 
\cite[Cor. 3.6]{Jardinejpaa}, \cite[\S 6.1]{Jardine}, and 
\cite[\S\S 4, 5]{hGal}, respectively. In several of these 
examples, a Postnikov tower provides an alternative to the 
hypercohomology spectrum that we use. In all of these 
constructions of the descent spectral sequence, some kind of 
finiteness assumption is required in order to identify the 
homotopy groups of the abutment as being those of the homotopy 
fixed point spectrum.
\end{Rk}
\par
Let $X$ be a discrete $G$-spectrum. We now develop a second 
model for $X^{hK}$, where $K$ is a closed subgroup of $G$, 
that is functorial in $K$.
\par
The map $X \rightarrow X_{f,G}$ in $Sp_K$ gives a weak 
equivalence $X^{hK} \rightarrow (X_{f,G})^{hK}.$ Composition 
with the weak equivalence $(X_{f,G})^{hK} \rightarrow 
\mathrm{holim} \, _\Delta \, 
(\Gamma^\bullet_K((X_{f,G})_{f,K}))^K$ gives a weak equivalence 
$X^{hK} \rightarrow \mathrm{holim} \, _\Delta \, 
(\Gamma^\bullet_K((X_{f,G})_{f,K}))^K$ between fibrant spectra. 
The inclusion $K \rightarrow G$ induces a morphism 
$\Gamma_G(X_{f,G}) \rightarrow \Gamma_K(X_{f,G}),$ giving a 
map $\Gamma^\bullet_G(X_{f,G}) \rightarrow 
\Gamma^\bullet_K(X_{f,G})$ of cosimplicial discrete $K$-spectra.
\begin{Lem}\label{zig-zag}
There is a weak equivalence \[\rho \: \mathrm{holim} \, 
_\Delta \, (\Gamma^\bullet_G(X_{f,G}))^K \negthinspace 
\rightarrow 
\mathrm{holim} \, _\Delta \, (\Gamma^\bullet_K(X_{f,G}))^K 
\negthinspace \rightarrow \mathrm{holim} \, _\Delta \, 
(\Gamma^\bullet_K((X_{f,G})_{f,K}))^K.\]
\end{Lem}  
\begin{proof}
Recall the conditionally convergent spectral sequence 
\[H^s_c(K;\pi_t(X)) \cong H^s_c(K; \pi_t((X_{f,G})_{f,K})) 
\Rightarrow \pi_{t-s}(\mathrm{holim} \, _\Delta \, 
(\Gamma^\bullet_K((X_{f,G})_{f,K}))^K).\]
\par
We compare this spectral sequence with the homotopy spectral 
sequence for $\mathrm{holim} \, _\Delta \, 
(\Gamma^\bullet_G(X_{f,G}))^K$. Note that if $Y$ is a 
discrete $G$-spectrum that is fibrant as a spectrum, then 
$\mathrm{Map}_c(G,Y) \cong \mathrm{colim} \, _{N} \, \prod 
\, _{G/N} \, Y$ and \[\mathrm{Map}_c(G,Y)^K \cong 
\mathrm{Map}_c(G/K,Y) \cong \mathrm{colim} \, _N \, 
\textstyle{\prod} \, _{G/(NK)} \, Y\] are fibrant spectra. 
Thus, $(\Gamma^\bullet_G(X_{f,G}))^K$ is a cosimplicial 
fibrant spectrum, and there is a conditionally convergent 
spectral sequence \[E_2^{s,t} = H^s((\Gamma^\ast_G 
\pi_t(X_{f,G}))^K) \Rightarrow \pi_{t-s}(\mathrm{holim} \, 
_\Delta \, (\Gamma^\bullet_G(X_{f,G}))^K).\] As in the proof 
of Theorem \ref{fibrantmodel}, $0 \rightarrow \pi_t(X_{f,G}) 
\rightarrow \Gamma^\ast_G(\pi_t(X_{f,G}))$ is a $(-)^K$-acyclic 
resolution of $\pi_t(X_{f,G})$. Hence, \[E_2^{s,t} \cong 
H^s_c(K;\pi_t(X_{f,G})) \cong H^s_c(K;\pi_t(X)).\]  
\par
Since $\rho$ is compatible with the isomorphism between the two 
$E_2$-terms, the spectral sequences are isomorphic and $\rho$ is 
a weak equivalence.
\end{proof}
\begin{Rk}\label{model2}
Lemma \ref{zig-zag} gives the following weak equivalences between 
fibrant spectra: \[X^{hK} = (X_{f,K})^K \rightarrow \mathrm{holim} 
\, _\Delta \, (\Gamma^\bullet_K((X_{f,G})_{f,K}))^K \leftarrow 
\mathrm{holim} \, _\Delta \, (\Gamma^\bullet_G(X_{f,G}))^K.\] 
Thus, 
if $\mathrm{vcd}(G) < \infty$, $X$ is a discrete $G$-spectrum, 
and $K$ is a closed subgroup of $G$, then 
$\mathrm{holim} \, _\Delta \, 
(\Gamma^\bullet_G(X_{f,G}))^K$ is a model for $X^{hK}$, so that 
\[X^{hK} =  \mathrm{holim} \, _\Delta \, 
(\Gamma^\bullet_G(X_{f,G}))^K\] is another definition of the 
homotopy fixed points.
\end{Rk}
\par
This discussion yields the following result.
\begin{Thm}\label{functor}
If $X$ is a discrete $G$-spectrum, with $\mathrm{vcd}(G) < 
\infty$, then there is a presheaf of spectra $P(X) \: 
(\mathcal{O}_G)^\mathrm{op} \rightarrow Sp$, defined by 
\[P(X)(G/K) = \mathrm{holim} \, _\Delta \, 
(\Gamma^\bullet_G(X_{f,G}))^K =X^{hK}.\]
\end{Thm} 
\begin{proof}
If $Y$ is a discrete $G$-set, any morphism $f \colon G/H 
\rightarrow G/K$, in $\mathcal{O}_G$, induces a map 
\[\mathrm{Map}_c(G,Y)^K \cong \mathrm{Map}_c(G/K,Y) \rightarrow 
\mathrm{Map}_c(G/H,Y) \cong \mathrm{Map}_c(G,Y)^H.\] Thus, if $Y 
\in Sp_G,$ $f$ induces a map $\mathrm{Map}_c(G,Y)^K \rightarrow 
\mathrm{Map}_c(G,Y)^H,$ so that there is a map \[P(X)(f) \: 
\mathrm{holim} 
\, _\Delta \, (\Gamma^\bullet_G(X_{f,G}))^K \rightarrow 
\mathrm{holim} \, _\Delta \, (\Gamma^\bullet_G(X_{f,G}))^H.\] 
It is easy to check that $P(X)$ is actually a functor.
\end{proof}
\par
We conclude this section by pointing out a useful fact: 
smashing with a finite spectrum, with trivial $G$-action, 
commutes with taking homotopy fixed points. To state this 
precisely, we first of all define the relevant map.
\par
Let $X$ be a discrete $G$-spectrum and let $Y$ be any spectrum 
with trivial $G$-action. Then there is a map \[(\mathrm{holim} 
\, _\Delta \, (\Gamma^\bullet_{G} X_f)^{G}) \wedge Y 
\rightarrow \mathrm{holim} \, _\Delta \, 
((\Gamma^\bullet_{G} X_f)^{G} \wedge Y) \rightarrow 
\mathrm{holim} \, _\Delta \, ((\Gamma^\bullet_{G} X_f) 
\wedge Y)^G.\] Also, there is a natural $G$-equivariant map 
$\mathrm{Map}_c(G,X) \wedge Y \rightarrow 
\mathrm{Map}_c(G,X \wedge Y)$ that is defined by the 
composition \[(\mathrm{colim} \, _N \, 
\textstyle{\prod} \, _{G/N} \, X) \wedge Y \cong 
\mathrm{colim} \, _N \, ((\textstyle{\prod} \, _{G/N} \, X) 
\wedge Y) \rightarrow \mathrm{colim} \, _N \, 
\textstyle{\prod} \, _{G/N} \, (X \wedge Y),\] by using the 
isomorphism $\mathrm{Map}_c(G,X) \cong \mathrm{colim} \, _N \, 
\textstyle{\prod} \, _{G/N} \, X.$ This gives a natural 
$G$-equivariant map \[(\Gamma_G \, \Gamma_G X) \wedge Y 
\rightarrow \Gamma_G \, ((\Gamma_G X) \wedge Y) \rightarrow 
\Gamma_G \, \Gamma_G (X \wedge Y).\] Thus, iteration gives a 
$G$-equivariant map of cosimplicial spectra 
\[(\Gamma^\bullet X) \wedge Y \rightarrow 
\Gamma^\bullet(X \wedge Y).\] Therefore, if 
$\mathrm{vcd}(G)< \infty,$ there is a canonical 
map $X^{hG}
\wedge Y \rightarrow (X_f \wedge Y)^{hG}$ that is 
defined by composing the map $X^{hG}
\wedge Y \rightarrow \mathrm{holim} \, _\Delta \, 
((\Gamma^\bullet X_f) \wedge Y)^G,$ from above, with 
the map \[\mathrm{holim} \, _\Delta \, ((\Gamma^\bullet 
X_f) \wedge Y)^G \rightarrow \mathrm{holim} \, _\Delta \, 
(\Gamma^\bullet (X_f \wedge Y))^G \rightarrow \mathrm{holim} 
\, _\Delta \, (\Gamma^\bullet (X_f \wedge Y)_f)^G.\]
\begin{Lem}[{\cite[Prop. 3.10]{Mitchell}}]\label{commute3}
If $\, \mathrm{vcd}(G) < \infty,$ $X \in Sp_G,$ and $Y$ is a
finite spectrum with trivial $G$-action, then $X^{hG}
\wedge Y \rightarrow (X_{f,G} \wedge Y)^{hG}$ is a weak 
equivalence.
\end{Lem}
\begin{Rk}\label{commute2}
By Lemma \ref{commute3}, when $Y$ is a finite spectrum, there 
is a zigzag of natural weak equivalences $X^{hG} \wedge Y 
\rightarrow (X_{f,G} \wedge Y)^{hG} \leftarrow (X \wedge 
Y)^{hG}.$ We refer to this zigzag by writing $X^{hG} \wedge 
Y \cong (X \wedge Y)^{hG}.$
\end{Rk}
\section{Homotopy fixed points for towers in $Sp_G$}\label{hfps2}
In this section, $\{Z_i\}$ is always a tower in $Sp_G$ (except in 
Definition \ref{M-L}). For $\{Z_i\}$ a tower of fibrant spectra, 
we define the homotopy fixed point spectrum $(\mathrm{holim} \, _i 
\, Z_i)^{hG}$ and construct its descent spectral sequence.
\begin{Def}\label{general}
If $\{Z_i\}$ in $\mathbf{tow}(Sp_G)$ is a tower of fibrant spectra, 
we 
define $Z = \mathrm{holim} \, _i \, Z_i$, a continuous
$G$-spectrum. The homotopy
fixed point spectrum $Z^{hG}$ is defined to be $\mathrm{holim} \, 
_i \, Z_i^{hG},$ a fibrant spectrum. 
\end{Def}
\par
We make several comments about Definition \ref{general}. Let $H$ 
be a closed subgroup of $G.$ Then the maps \[\mathrm{holim} \,
_i \, ((Z_i)_f)^H \rightarrow \mathrm{holim} \, _i \,
\mathrm{holim} \, _\Delta \, (\Gamma^\bullet_H(Z_i)_f)^H \ \ \ 
\mathrm{and}\] \[\mathrm{holim} \, _i \,
\mathrm{holim} \, _\Delta \, (\Gamma^\bullet_G(Z_i)_{f,G})^H 
\rightarrow
\mathrm{holim} \, _i \,
\mathrm{holim} \, _\Delta \, 
(\Gamma^\bullet_H((Z_i)_{f,G})_{f,H})^H\] are weak equivalences. 
Thus, in Definition \ref{general}, any one of our three 
definitions 
for homotopy fixed points (Definition \ref{profinite}, Remarks 
\ref{model} and \ref{model2}) can be used for $Z_i^{hH}$.
\par
In Definition \ref{general}, suppose that not all the $Z_i$ are 
fibrant in $Sp$. Then the map $Z=\mathrm{holim} \, _i \, Z_i 
\rightarrow \mathrm{holim} \, _i \, (Z_i)_{f,\{e\}} = 
\mathrm{holim} \, _i \, Z_i^{h\{e\}} = Z^{h\{e\}}$ need not be 
a weak equivalence. Thus, for an arbitrary tower in $Sp_G$, 
Definition \ref{general} can fail to have the desired property 
that $Z \rightarrow Z^{h\{e\}}$ is a weak equivalence. 
\par
Below, Lemmas \ref{special} and \ref{ex1}, and Remark 
\ref{ex1remark}, show that when $G$ is a finite group, $Z^{hG} 
\simeq Z^{h'G}$, and, for any $G$, $Z^{hG}$ can be obtained by 
using a total right derived functor that comes from fixed points. 
Thus, Definition \ref{general} generalizes the notion of homotopy 
fixed points to towers of discrete $G$-spectra.
\begin{Lem}\label{special}
Let $G$ be a finite group and let $\{Z_i\}$ in 
$\mathbf{tow}(Sp_G)$ be a tower of fibrant spectra. Then there 
is a weak equivalence $Z^{hG} \rightarrow Z^{h'G}.$
\end{Lem}
\begin{proof}
We have: $Z^{hG} = \mathrm{holim} \, _i \, \mathrm{lim} \, _G \, 
(Z_i)_f$ and $Z^{h'G} = \mathrm{holim} \, _G \, \mathrm{holim} 
\, _i \,
(Z_i)_f$ (since $Z = \mathrm{holim} \, _i \, Z_i
\rightarrow \mathrm{holim} \, _i \, (Z_i)_f$ is a weak
equivalence and $G$-equivariant, and $\mathrm{holim} \, _i
\, (Z_i)_f$ is fibrant in $Sp$). Then the map $Z^{hG} 
\rightarrow 
Z^{h'G}$ is defined to be \[\mathrm{holim} \, _i \, 
\mathrm{lim} \, _G \, (Z_i)_f \rightarrow \mathrm{holim} \, 
_i \, \mathrm{holim} \, _G \,
(Z_i)_f \cong \mathrm{holim} \, _G \, \mathrm{holim} \, _i \,
(Z_i)_f.\]
For each $i$, $\mathrm{lim} \, _G \,
(Z_i)_f \rightarrow \mathrm{holim} \, _G \, (Z_i)_f$ is a weak
equivalence between fibrant spectra, so that the desired map 
is a
weak equivalence. 
\end{proof}
\par
In the lemma below, \[\mathbf{R}(\mathrm{lim} \, _i \, (-)^G) 
\: \mathrm{Ho}(\mathbf{tow}(Sp_G)) \rightarrow \mathrm{Ho}(Sp)\] 
is the total right derived functor of the functor $\mathrm{lim} 
\, _i \, (-)^G \: \mathbf{tow}(Sp_G) \rightarrow Sp.$
\begin{Lem}\label{ex1}
If $\{Z_i\}$ is an arbitrary tower in $Sp_G$, then 
\[\mathrm{holim} \, _i \, ((Z_i)_f)^G
  \overset{\simeq}{\longrightarrow} \mathrm{holim} \, _i \,
  ({(Z_i)'_{\negthinspace f}})^G \overset{\simeq}{\longleftarrow}
  \mathrm{lim} \, _i \, ({(Z_i)'_{\negthinspace f}})^G =
  \mathbf{R}(\mathrm{lim} \, _i \, (-)^G)(\{Z_i\}).\]
\end{Lem}
\begin{proof}
The first weak equivalence follows from $((Z_i)_f)^G \rightarrow 
((Z_i)'_{\negthinspace f})^G$ being a weak equivalence between 
fibrant spectra. The second weak
equivalence holds because $\{({(Z_i)'_{\negthinspace f}})^G\}$ 
is a tower of fibrant spectra, such that all maps in the tower 
are fibrations. Finally, the equality is because $\{Z_i\} 
\rightarrow
\{(Z_i)'_{\negthinspace f}\}$ is a trivial cofibration in 
$\mathbf{tow}(Sp_G)$. 
\end{proof}
\begin{Rk}\label{ex1remark}
By Lemma \ref{tright}, if $X \in Sp_G$, then 
$X^{hG} = (\mathbf{R}(-)^G)(X)$. Also, by Lemma \ref{ex1}, if 
$\{Z_i\}$ in $\mathbf{tow}(Sp_G)$ is a tower of fibrant spectra, 
then \[Z^{hG} = \mathrm{holim} \, _i \, Z_i^{hG} = 
\mathrm{holim} \, _i \, ((Z_i)_f)^G \cong \mathbf{R}(\mathrm{lim} 
\, _i \, (-)^G)(\{Z_i\}).\] Thus, the homotopy fixed point 
spectrum $Z^{hG}$ is again given by the total right derived 
functor of an appropriately defined functor involving $G$-fixed 
points.
\end{Rk}
\par
Given any tower in $Sp_G$ of fibrant spectra, there is a descent 
spectral sequence whose $E_2$-term is a version of continuous 
cohomology, whose definition we now recall from \cite{Jannsen}. 
We use $\mathrm{DMod}(G)$ to denote the category of discrete 
$G$-modules. 
\begin{Def}\label{Jannsen}
Let $\mathrm{DMod}(G)^\Nset$ denote the category of
diagrams in discrete $G$-modules of the form $ \cdots
\rightarrow M_2 \rightarrow M_1 \rightarrow M_0$. Then 
$H^s_\mathrm{cont}(G; \{M_i\})$, the continuous cohomology of 
$G$
with coefficients in the tower $\{M_i\}$, is the $s$th right 
derived functor of
the left exact functor \[\mathrm{lim} \, _i \, (-)^G \: 
\mathrm{DMod}(G)^\mathbb{N} \rightarrow \mathbf{Ab}, \ \ \
\{M_i\} \mapsto \mathrm{lim} \, _i \, M_i^G.\] By 
\cite[Theorem 2.2]{Jannsen}, if the tower of abelian groups 
$\{M_i\}$ satisfies the Mittag-Leffler condition, then 
$H^s_\mathrm{cont}(G;\{M_i\}) \cong 
H^s_\mathrm{cts}(G;\mathrm{lim} \, _i \, M_i),$ for $s \geq 0.$
\end{Def}
\begin{Thm}\label{generaldss}
If $\mathrm{vcd}(G)<\infty$ and $\{Z_i\}$ in $\mathbf{tow}(Sp_G)$ 
is a tower of fibrant spectra, then there is a conditionally 
convergent descent spectral sequence 
\begin{equation}\label{geisser2}\zig 
H^s_\mathrm{cont}(G;\{\pi_t(Z_i)\}) \Rightarrow 
\pi_{t-s}(Z^{hG}).\end{equation}
\end{Thm}
\begin{Rk}\label{l-adic}
Such a spectral sequence goes back to Thomason's construction 
of the $\ell$-adic descent spectral sequence of
algebraic $K$-theory (\cite{Thomason},
\cite{Mitchell}). Spectral sequence (\ref{geisser2}) is the 
homotopy spectral sequence \[E_2^{s,t} = 
\textstyle{\mathrm{lim}} ^s \, _{\Delta \times \{i\}} \, 
\pi_t((\Gamma^\bullet_G((Z_i)_{f,G}))^G) \Rightarrow 
\pi_{t-s}(\mathrm{holim} \, _{\Delta \times \{i\}} \, 
(\Gamma^\bullet_G((Z_i)_{f,G}))^G),\] where the identification 
$E_2^{s,t} \cong H^s_\mathrm{cont}(G;\{\pi_t(Z_i)\})$ depends 
on the fact that if $\{I_i\}$ is an injective object in 
$\mathrm{DMod}(G)^\mathbb{N}$, then $\textstyle{\mathrm{lim}} 
^s \, _{\Delta \times \{i\}} \, (\Gamma^\bullet_G I_i)^G =0$, 
for $s >0$. We omit the details of the proof, since it is a 
special case of \cite[Prop. 3.1.2]{Geisser}, and also because 
(\ref{geisser2}) is not our focus of interest.
\end{Rk}
\par
For our applications, instead of spectral sequence 
(\ref{geisser2}), we are more interested in descent spectral 
sequence (\ref{geisser}) below. Spectral sequence (\ref{geisser}), 
a homotopy spectral sequence for a particular cosimplicial 
spectrum, 
is more suitable for comparison with the $K(n)$-local $E_n$-Adams 
spectral sequence (see \cite[Prop. A.5]{DH}), when (\ref{geisser}) 
has abutment $\pi_\ast((E_n \wedge X)^{hG})$, where $X$ is a 
finite spectrum.
\begin{Def}\label{M-L}
If $\{Z_i\}$ is a tower of spectra such that $\{\pi_t(Z_i)\}$ 
satisfies the Mittag-Leffler condition for every $t \in 
\Zset$, then $\{Z_i\}$ is a {\em Mittag-Leffler tower\/} of 
spectra.
\end{Def}
\begin{Thm}\label{Dss}
If $\mathrm{vcd}(G) < \infty$ and $\{Z_i\}$ in 
$\mathbf{tow}(Sp_G)$ is a tower of fibrant spectra, then
there is a conditionally convergent descent spectral sequence
\begin{equation}\label{geisser}\zig
E_2^{s,t} = \pi^s\pi_t (\mathrm{holim} \,
_i \, (\Gamma^\bullet_G (Z_i)_f)^G)  \Rightarrow
\pi_{t-s}(Z^{hG}).\end{equation} If $\{Z_i\}$ is a 
Mittag-Leffler tower, then $E_2^{s,t} \cong 
H^s_\mathrm{cts}(G; \pi_t(Z))$. 
\end{Thm}
\begin{Rk}\label{differ}
In Theorem \ref{Dss}, when $\{Z_i\}$ is a Mittag-Leffler tower, 
$E_2^{s,t} \cong H^s_\mathrm{cts}(G;\pi_t(Z)) \cong 
H^s_\mathrm{cont}(G;\{\pi_t(Z_i)\}),$ and spectral sequence 
(\ref{geisser}) is identical to (\ref{geisser2}). However, in 
general, spectral sequences (\ref{geisser2}) and (\ref{geisser}) 
are different. For example, if $G=\{e\}$, then in 
(\ref{geisser2}), $E_2^{0,t}=\mathrm{lim} \, _i \, \pi_t(Z_i)$, 
whereas in (\ref{geisser}), $E_2^{0,t}=\pi_t(\mathrm{holim} \, 
_i \, Z_i)$. 
\end{Rk}
\begin{proof}[Proof of Theorem \ref{Dss}.]
Note that $Z^{hG} \cong \mathrm{holim} \, _\Delta \, 
\mathrm{holim} \,
_i \, (\Gamma^\bullet_G (Z_i)_f)^G,$ and the diagram $\mathrm{holim} \, _i \, 
(\Gamma^\bullet_G (Z_i)_f)^G$ is a cosimplicial fibrant spectrum.
\par
Let $\{Z_i\}$ be a Mittag-Leffler tower. For $k \geq 0$, Lemma 
\ref{mittag} and Remark \ref{mittagrk} imply that 
$\textstyle{\mathrm{lim}}^1 \, _i \,
  \mathrm{Map}_c(G^{k}, \pi_{t+1}(Z_i))  = 0.$ Therefore,
\[\pi_t (\mathrm{holim} \, _i \,
  (\mathrm{Map}_c(G^{k+1}, (Z_i)_f ))^G ) \cong \mathrm{lim} \, 
_i
  \, \mathrm{Map}_c(G^{k+1}, \pi_t(Z_i))^G,\] and hence, 
$\pi_t(\mathrm{holim} \, _i \,
  (\Gamma^\bullet (Z_i)_f)^G) \cong \mathrm{lim} \, _i \,
  (\Gamma^\bullet \pi_t(Z_i))^G.$ Thus, \[E_2^{s,t} \cong \pi^s 
(\mathrm{lim} \, _i \,
  (\Gamma^\bullet \pi_t(Z_i))^G) \cong H^s(\mathrm{lim} \, _i \, 
(-)^G \, \{\Gamma^\ast \pi_t(Z_i)\}_i).\] 
\par
Consider the exact sequence $\{0\} \rightarrow \{\pi_t(Z_i)\} 
\rightarrow \{\Gamma^\ast\pi_t(Z_i)\}$ in 
$\mathrm{DMod}(G)^\mathbb{N}.$ Note that, for $s,k>0$, by Theorem 
\ref{ses}, \[H^s_\mathrm{cont}(G;\{\Gamma^k \pi_t(Z_i)\}) \cong 
H^s_\mathrm{cts}(G;\mathrm{lim} \, _i \, \Gamma^k \pi_t(Z_i)) 
\cong \mathrm{lim} \, _i \, H^s_c(G; \Gamma^k \pi_t(Z_i))=0,\] 
since the tower $\{\Gamma^k \pi_t(Z_i)\}$ satisfies the 
Mittag-Leffler condition, and, for each $i$, $\Gamma^k \pi_t(Z_i) 
\cong \mathrm{Map}_c(G,\Gamma^{k-1} \pi_t(Z_i))$ is 
$(-)^G$-acyclic. Thus, the above exact sequence is a $(\mathrm{lim} 
\, _i \, (-)^G)$-acyclic resolution of $\{\pi_t(Z_i)\},$ so that 
we obtain the isomorphism $E_2^{s,t} \cong 
H^s_\mathrm{cont}(G;\{\pi_t(Z_i)\}) \cong H^s_\mathrm{cts}(G; 
\pi_t(Z)).$
\end{proof}
\par
By Remark \ref{ex1remark}, we can rewrite spectral sequence 
(\ref{geisser}), when $\{Z_i\}$ is a Mittag-Leffler tower, in 
a more
conceptual way: \[R^s(\mathrm{lim} \, _i \,
(-)^G)\{\pi_t(Z_i)\} \Rightarrow
\pi_{t-s}(\mathbf{R}(\mathrm{lim} \, _i \,
(-)^G)(\{Z_i\})).\] Spectral sequence (\ref{geisser2}) can always 
be written in this way.
\section{Homotopy fixed point spectra for $\Lhat  (E_n \wedge X)$}
\par
In this section, for any spectrum $X$ and any closed subgroup $G$ 
of $G_n$, we define the homotopy fixed point spectrum 
$(\Lhat(E_n \wedge X))^{hG}$, using the continuous action of $G$.
\par
Let $X$ be an arbitrary spectrum with trivial $G_n$-action. By 
Corollary \ref{Rezk}, there is a weak equivalence $F_n \wedge M_I 
\wedge X \rightarrow E_n \wedge M_I \wedge X.$ Then, using 
functorial fibrant replacement, there is a map of towers 
\[\{(F_n \wedge M_I \wedge 
X)_\mathtt{f}\} \rightarrow \{(E_n \wedge M_I \wedge 
X)_\mathtt{f}\},\] which yields the weak equivalence \[\Lhat 
(E_n \wedge X) \cong \mathrm{holim} \, _I \, (E_n \wedge M_I 
\wedge X)_\mathtt{f} \overset{\simeq}{\longleftarrow} 
\mathrm{holim} \, _I \, (F_n \wedge M_I \wedge X)_\mathtt{f}.\] 
As in the proof of Theorem \ref{wonder}, this implies the 
following lemma, since the diagram $\{(F_n \wedge M_I 
\wedge X)_{f,G_n}\}$ is a tower of 
fibrant spectra.
\begin{Lem}\label{gencts}
Given any spectrum $X$ with trivial $G_n$-action, the isomorphism 
\[\Lhat (E_n \wedge X) \cong \mathrm{holim} \, _I \, (F_n \wedge 
M_I \wedge X)_{f,G_n}\] makes $\Lhat(E_n \wedge X)$ a continuous 
$G_n$-spectrum.
\end{Lem}
Let $G$ be any closed subgroup of $G_n$. Since $\{(F_n \wedge 
M_I \wedge X)_{f,G_n}\}$ is a tower of discrete $G$-spectra that 
are fibrant in $Sp$, Lemma \ref{gencts} also shows that 
$\Lhat(E_n \wedge X)$ is a continuous $G$-spectrum. By Theorem 
\ref{finally}, the composition \[(F_n \wedge M_I \wedge X) 
\rightarrow  (F_n \wedge M_I \wedge X)_{f,G_n} \rightarrow 
((F_n \wedge M_I \wedge X)_{f,G_n})_{f,G}\] is a trivial 
cofibration in $Sp_G$, with target fibrant in $Sp_G$. Therefore, 
\[((F_n \wedge M_I \wedge X)_{f,G_n})^{hG} = (F_n \wedge M_I 
\wedge X)^{hG}.\] 
Then, by Definition \ref{general}, we have the following.
\begin{Def}\label{bigdef}
Let $G<_c G_n$. Then \[E_n^{hG} = (\mathrm{holim} \, _I \, 
(F_n \wedge M_I)_{f,G_n})^{hG} = \mathrm{holim} \, _I \, (F_n 
\wedge M_I)^{hG}.\] More generally, for any spectrum $X$, 
\[(\Lhat (E_n \wedge X))^{hG} = (\mathrm{holim} \, _I \, 
(F_n \wedge M_I \wedge X)_{f,G_n})^{hG} = \mathrm{holim} \, 
_I \, (F_n \wedge M_I \wedge X)^{hG}.\]
\end{Def}
\begin{Rk}
When $X$ is a finite spectrum, $E_n \wedge X \simeq \Lhat 
(E_n \wedge X)$. Thus, we have $(E_n \wedge X)^{hG} \cong 
(\Lhat (E_n \wedge X))^{hG}.$
\end{Rk}
\begin{Rk}\label{elementary}
For any $X$, $\Lhat(E_n \wedge X) \cong \mathrm{holim} \, _I 
\, (F_n \wedge M_I \wedge X)_{f,G}$ also shows that $\Lhat(E_n 
\wedge X)$ is a continuous $G$-spectrum. Note that, by 
definition, the spectra $((F_n \wedge M_I \wedge X)_{f,G})^{hG}$ 
and $(F_n \wedge M_I \wedge X)^{hG}$ are identical. This implies 
that \[(\Lhat(E_n \wedge X))^{hG} = \mathrm{holim} \, _I \, 
((F_n \wedge M_I \wedge X)_{f,G})^{hG} =\mathrm{holim} \, _I 
\, (F_n \wedge M_I \wedge X)^{hG},\] as before.
\end{Rk}
\par
Note that Definition \ref{bigdef} implies the identifications 
\[E_n^{hG} = \mathrm{holim} \,
_I \, \mathrm{holim} \, _\Delta \, (\Gamma^\bullet_{G_n} (F_n 
\wedge
M_I)_{f,G_n})^G\] and
\[(\Lhat (E_n \wedge X))^{hG} = \mathrm{holim} \,
_I \, \mathrm{holim} \, _\Delta \, (\Gamma^\bullet_{G_n} (F_n 
\wedge
M_I \wedge X)_{f,G_n})^G.\]
The first identification, coupled with Theorem \ref{functor}, 
implies the following result.
\begin{Thm}
There is a functor \[P \: (\mathcal{O}_{G_n})^\mathrm{op} 
\rightarrow Sp, \ \ \ P(G_n/G) = E_n^{hG},\] where $G$ is any 
closed subgroup of $G_n$.
\end{Thm}  
\par
In addition to the above identifications, we also have 
\[E_n^{hG} = \mathrm{holim} \,
_I \, \mathrm{holim} \, _\Delta \, (\Gamma^\bullet_{G} (F_n 
\wedge
M_I)_{f,G})^G\] and
\[(\Lhat (E_n \wedge X))^{hG} = \mathrm{holim} \,
_I \, \mathrm{holim} \, _\Delta \, (\Gamma^\bullet_{G} (F_n 
\wedge
M_I \wedge X)_{f,G})^G.\] 
\par
Since $E_n^{\mathtt{h}G}$ is $K(n)$-local, one would expect that 
$E_n^{hG}$ is $K(n)$-local; this is verified below.
\begin{Lem}\label{local}
Let $G$ be a closed subgroup of $G_n$ and let $X$ be any spectrum. 
Then the spectra $(F_n \wedge M_I \wedge 
X)^{hG}$ and
$(\Lhat(E_n \wedge X))^{hG}$ are $K(n)$-local. Also, 
$(F_n \wedge X)^{hG}$ is 
$E(n)$-local.
\end{Lem}
\begin{proof}
Recall that $(F_n \wedge M_I \wedge X)^{hG} = \mathrm{holim} \, 
_\Delta \, (\Gamma^\bullet_G (F_n \wedge M_I \wedge X)_{f,G})^G,$ 
and, for $k \geq 1,$  
\[(\Gamma^k_G(F_n \wedge M_I \wedge X)_f)^G \simeq 
\mathrm{Map}_c(G^{k-1}, F_n \wedge X) \wedge M_I.\] Since 
$\mathrm{Map}_c(G, F_n \wedge X) \cong \mathrm{colim} \, _i
\, \textstyle{\prod} \, _{G/(U_i \cap G)} \, (F_n \wedge X)$, 
and $F_n \wedge X$ is $E(n)$-local, the finite product is too, 
and hence, the direct limit $\mathrm{Map}_c(G, F_n \wedge X)$ 
is $E(n)$-local. Iterating this argument shows that 
$\mathrm{Map}_c(G^{k-1}, F_n \wedge X) \cong \Gamma_G 
\Gamma_G \cdots \Gamma_G (F_n \wedge X)$ is $E(n)$-local. 
Smashing $\mathrm{Map}_c(G^{k-1}, F_n \wedge X)$ with the 
type $n$ spectrum $M_I$ shows that $(\Gamma^k_G(F_n \wedge 
M_I \wedge X)_f)^G$
is $K(n)$-local. Therefore, since the homotopy limit of an 
arbitrary 
diagram of $E$-local
spectra is $E$-local \cite[pg. 259]{Bousfieldlocal}, 
$(F_n \wedge M_I \wedge X)^{hG}$ and $(\Lhat(E_n \wedge X))^{hG}$ 
are $K(n)$-local. The same argument shows that 
$(F_n \wedge X)^{hG}$ is 
$E(n)$-local. 
\end{proof}
\par
The following theorem shows that the homotopy fixed points of 
$\Lhat(E_n \wedge X)$ are obtained by taking the 
$K(n)$-localization of the homotopy fixed points of the 
discrete $G_n$-spectrum 
$(F_n \wedge X)$. Thus, homotopy fixed points for $E_n$ come from 
homotopy fixed points for $F_n$.
\begin{Thm}\label{fun}
For any closed subgroup $G$ of $G_n$ and any spectrum $X$ with 
trivial $G$-action, there is an isomorphism \[(\Lhat(E_n \wedge 
X))^{hG} \cong \Lhat ((F_n \wedge X)^{hG})\] in the stable 
category. In particular, $E_n^{hG} \cong \Lhat(F_n^{hG})$.
\end{Thm}
\begin{proof}
After switching $M_I$ and $X$, $(\Lhat(E_n \wedge X))^{hG} 
\cong \mathrm{holim} \, 
_I \, (F_n \wedge X \wedge M_I)^{hG}$. By Remark \ref{commute2}, 
$(F_n \wedge X \wedge M_I)^{hG} \cong (F_n \wedge X)^{hG} 
\wedge M_I \simeq ((F_n \wedge X)^{hG} \wedge M_I)_\mathtt{f},$ 
where the isomorphism signifies a zigzag of natural weak 
equivalences. Thus, \[(\Lhat(E_n \wedge X))^{hG} \cong 
\mathrm{holim} \, _I \, 
((F_n \wedge X)^{hG} \wedge M_I)_\mathtt{f} \cong \Lhat((F_n 
\wedge X)^{hG}),\] since, by Lemma \ref{local}, 
$(F_n \wedge X)^{hG}$ is $E(n)$-local.
\end{proof}
\begin{Cor}\label{different} If $X$ is a finite spectrum
  of type $n$, then there is an isomorphism $(F_n
  \wedge X)^{hG} \cong E_n^{hG}
  \wedge X.$ In particular, $(F_n \wedge M_I)^{hG} \cong E_n^{hG} 
\wedge M_I.$\end{Cor}
\begin{proof}
Since $X$ is finite, by Remark \ref{commute2}, $(F_n
  \wedge X)^{hG} \cong F_n^{hG} \wedge X.$ Then, since $F_n^{hG}$ 
is $E(n)$-local by Lemma \ref{local}, there are isomorphisms 
$F_n^{hG} \wedge X \cong \Lhat (F_n^{hG}) \wedge X 
\cong E_n^{hG} \wedge X.$
\end{proof} 
\par
We conclude this section by observing that smashing with a 
finite spectrum 
commutes with taking the homotopy fixed points of $E_n$. 
\begin{Thm}\label{final}
Let $G$ be a closed subgroup of $G_n$ and let $X$ be a finite 
spectrum. 
Then there is an isomorphism $(E_n \wedge X)^{hG} \cong 
E_n^{hG} \wedge X.$
\end{Thm}
\begin{proof}
Recall that $(E_n \wedge X)^{hG} = \mathrm{holim} \, _I \, 
(F_n \wedge M_I \wedge X)^{hG}.$ The zigzag of natural weak 
equivalences between $(F_n \wedge M_I \wedge X)^{hG}$ and $(F_n 
\wedge M_I)^{hG} \wedge X$ yields \[(E_n \wedge X)^{hG} 
\cong \mathrm{holim} \, _I \, ((F_n \wedge M_I)^{hG} \wedge 
X)_\mathtt{f} \simeq (\mathrm{holim} \, _I \, (F_n \wedge 
M_I)^{hG}) \wedge X,\] where the weak equivalence is due to 
Lemma \ref{commute}.
\end{proof}
\section{The descent spectral sequence for $(\Lhat(E_n \wedge 
X))^{hG}$}
\par
By applying the preceding two sections, it is now an easy matter 
to build the descent spectral sequence for $(\Lhat(E_n \wedge 
X))^{hG}$.
\begin{Def}
Let $X$ be a spectrum. If the tower $\{\pi_t(E_n \wedge M_I 
\wedge X)\}_I$ of abelian groups satisfies the 
Mittag-Leffler condition for all $t \in \Zset$, then $X$ 
is an {\em $E_n$-Mittag-Leffler spectrum}. If 
$X$ is an $E_n$-Mittag-Leffler spectrum, then, for convenience, 
we say that 
$X$ is $E_n$-ML.
\end{Def}
\par
Any finite spectrum $X$ is $E_n$-ML, since $\{\pi_t(E_n \wedge 
M_I \wedge X)\}_I$ is a tower of finite
abelian groups, by Lemma \ref{finite}. However, an
$E_n$-ML spectrum need not be finite. For example, for $j \geq
1$, let $X=E_n^{(j)}$. Then $\pi_t(E_n \wedge M_I
\wedge X) \cong
\mathrm{Map}_c^\ell(G_n^{j}, \pi_t(E_n)/I).$ Since
$\{\pi_t(E_n)/I\}$ is a tower of epimorphisms, the tower 
$\{\mathrm{Map}_c^\ell(G_n^{j}, \pi_t(E_n)/I)\}$ is also, and 
$E_n^{(j)}$ is $E_n$-ML. 
\begin{Thm}\label{label}
Let $G$ be a closed subgroup of $G_n$ and let $X$ be any spectrum 
with trivial $G$-action. Let \[E_2^{s,t}= \pi^s\pi_t(\mathrm{holim} 
\, _I \, (\Gamma^\bullet_{G} (F_n \wedge
M_I \wedge X)_{f,G})^G).\] Then there is a conditionally convergent 
descent spectral sequence
\begin{equation}\zig\label{mydss}
E_2^{s,t} \Rightarrow \pi_{t-s}((\Lhat (E_n \wedge X) )^{hG}).
\end{equation}
If $X$ is $E_n$-ML, then 
$E_2^{s,t}\cong H^s_\mathrm{cts}(G;\pi_t(\Lhat(E_n \wedge X))).$
In particular, if $X$ is a finite spectrum, then descent 
spectral
sequence (\ref{mydss}) has the form
\[H^s_c(G; \pi_t(E_n 
\wedge X)) \Rightarrow
\pi_{t-s}((E_n \wedge X)^{hG}).\]
\end{Thm}
\begin{proof}
As in Remark \ref{elementary}, $\Lhat(E_n \wedge X) \cong 
\mathrm{holim} \, _I \, (F_n 
\wedge M_I \wedge X)_{f,G}$ is a continuous $G$-spectrum. Then 
(\ref{mydss}) follows from Theorem \ref{Dss} by considering the 
tower of spectra
$\{(F_n \wedge M_I \wedge X)_f \}_I.$ When $X$ is $E_n$-ML, 
$\{(F_n \wedge M_I \wedge X)_f \}$ is a Mittag-Leffler tower of 
spectra, and thus, the simplification of the $E_2$-term in this 
case follows from Theorem \ref{Dss}. By Corollary
\ref{help2}, when $X$ is finite, there is an isomorphism
$H^s_\mathrm{cts}(G;\pi_t(E_n \wedge X)) \cong H^s_{c}(G;
\pi_t(E_n \wedge X)).$ \end{proof}
\par
As discussed in Remark \ref{differ}, Theorem \ref{generaldss} 
gives a spectral sequence with abutment $\pi_\ast((\Lhat(E_n 
\wedge X))^{hG})$, the same as the abutment of (\ref{mydss}), but 
with $E_2$-term $H^s_\mathrm{cont}(G;\{\pi_t(E_n \wedge M_I \wedge 
X)\}),$ which is, in general, different from the $E_2$-term of 
(\ref{mydss}). We are interested in the descent spectral sequence 
of Theorem \ref{label}, not just because it is a second descent 
spectral sequence with an interesting $E_2$-term, but, as 
mentioned 
in \S \ref{hfps2}, it can be compared with the $K(n)$-local 
$E_n$-Adams spectral sequence. (This comparison is work in 
progress.)
\par
We conclude this paper with a computation that uses spectral 
sequence (\ref{mydss}) to compute $\pi_\ast((\Lhat(E_n \wedge 
E_n^{(j)}))^{hG_n}),$ where $j \geq 1$ and $G_n$ acts only on the
leftmost factor. By Theorem \ref{powers}, 
\begin{align*}\pi_t(\Lhat  (E_n^{(j+1)})) & \cong 
\mathrm{Map}_c^\ell(G_n^j, \pi_t(E_n)) \\ & \cong 
\mathrm{lim} \, _I \, \mathrm{Map}_c^\ell(G_n, 
\mathrm{Map}_c(G_n^{j-1}, \pi_t(E_n \wedge M_I))) \\ 
& \cong \mathrm{lim} \, _I \, \mathrm{Map}_c(G_n, 
\mathrm{Map}_c(G_n^{j-1}, \pi_t(E_n \wedge M_I))),\end{align*} 
where the tower $\{\mathrm{Map}_c(G_n, \mathrm{Map}_c(G_n^{j-1}, 
\pi_t(E_n \wedge M_I)))\}$ satisfies the Mittag-Leffler condition, 
by Remark \ref{mittagrk}. 
Thus, in spectral sequence (\ref{mydss}), Theorem \ref{ses} 
implies 
that
\begin{align*}E_2^{s,t} & \cong H^s_\mathrm{cts}(G_n; \pi_t(\Lhat  
(E_n^{(j+1)}))) \\ & \cong \mathrm{lim} \, _I \, H^s_c(G_n; 
\mathrm{Map}_c(G_n, \mathrm{Map}_c(G_n^{j-1}, \pi_t(E_n \wedge 
M_I)))),\end{align*} which vanishes for $s>0$, and
  equals $\mathrm{Map}_c(G_n^{j-1}, \pi_t(E_n)),$ when $s=0.$ 
Thus, \[\pi_\ast ((\Lhat  (E_n^{(j+1)}))^{hG_n}) \cong 
\mathrm{Map}_c(G_n^{j-1},
  \pi_\ast(E_n)) \cong \pi_\ast(\Lhat(E_n^{(j)})),\] as 
abelian groups. Therefore, for $j \geq 1$, there is an 
isomorphism 
$(\Lhat  (E_n^{(j+1)}))^{hG_n} \cong \Lhat  (E_n^{(j)}),$ 
in the stable category. The 
techniques described in this paper do not allow us to handle 
the $j=0$ case, which would say that $E_n^{hG_n}$ and 
$\Lhat(S^0) \simeq E_n^{\mathtt{h}G_n}$ are isomorphic.

\bibliographystyle{plain}
\bibliography{biblio}

\end{document}